\documentclass[11pt]{article}

\usepackage{latexsym,amsmath,amsfonts,amssymb,amstext,theorem,euscript,array,enumerate,mathrsfs}
\usepackage{color}
\usepackage{comment}

\newtheorem{Theorem}{Theorem}[section]
\newtheorem{Definition}{Definition}[section]
\newtheorem{Proposition}{Proposition}[section]

\newtheorem{Lemma}{Lemma}[section]

\newtheorem{Remark}{Remark}[section]

\numberwithin{equation}{section}

\def\esssup_#1{\underset{#1}{\mathrm{ess\,sup\, }}}
\def\essinf_#1{\underset{#1}{\mathrm{ess\,inf\, }}}

\def\qed{{\hfill\hbox{\enspace${ \square}$}} \smallskip}
\def\sqr#1#2{{\vcenter{\vbox{\hrule height .#2pt \hbox{\vrule
 width .#2pt height#1pt \kern#1pt \vrule
width .#2pt} \hrule height .#2pt}}}}
\def\square{\mathchoice\sqr54\sqr54\sqr{4.1}3\sqr{3.5}3}

\def\ds{\begin{displaystyle}}
\def\eds{\end{displaystyle}}
\def\dis{\displaystyle }
\def\<{\langle }
\def\>{\rangle }

\def \N{\mathbb{N}}
\def \R{\mathbb{R}}

\def \E{\mathbb{E}}
\def \F{\mathbb{F}}
\def \G{\mathbb{G}}
\def \P{\mathbb{P}}
\def \Q{\mathbb{Q}}

\def \Ac{{\cal A}}

\def \Ec{{\cal E}}
\def \Fc{{\cal F}}

\def \Kc{{\cal K}}

\def \Rc{{\cal R}}
\def \Sc{{\cal S}}
\def \Tc{{\cal T}}

\def \Vc{{\cal V}}

\def \Vc{{\cal V}}

\def\cala{{\cal A}}
\def\calb{{\cal B}}
\def\calc{{\cal C}}

\def\calf{{\cal F}}
\def\calg{{\cal G}}

\def\caln{{\cal N}}

\def\calp{{\cal P}}

\def\calv{{\cal V}}
\def\call{{\cal L}}

\def\bfB{{\bf B}}
\def\bfC{{\bf C}}

\def \eps{\varepsilon}

\def \ep{\hbox{ }\hfill$\Box$}

\def\reff#1{{\rm(\ref{#1})}}

\def\beqs{\begin{eqnarray*}}
\def\enqs{\end{eqnarray*}}
\def\beq{\begin{eqnarray}}
\def\enq{\end{eqnarray}}

\allowdisplaybreaks

\begin{document}

\title{Optimal switching problems
 with an infinite set of modes:
 \\
 an approach by randomization and constrained backward SDEs }

\author{
Marco FUHRMAN
\thanks{Universit\`{a} degli Studi di Milano, Dipartimento di Matematica, via Saldini 50, 20133 Milano, Italy; e-mail: \texttt{marco.fuhrman@unimi.it}}
\and
Marie-Am\'{e}lie MORLAIS
\thanks{Le Mans Universit\'{e}, Laboratoire Manceau de Math\'ematiques, Avenue Olivier Messiaen,
72085 Le Mans, Cedex 9, France; e-mail: \texttt{Marie\_Amelie.Morlais@univ-lemans.fr}
}
}

\maketitle

\begin{abstract}
We address a general optimal switching problem over finite horizon for a stochastic system described
by a differential equation driven by Brownian motion. The main novelty is the
fact that we allow for infinitely many modes (or regimes, i.e. the possible
values of the piecewise-constant control process).  We allow all the given coefficients in the model
to be path-dependent, that is, their value at any time depends on the past
trajectory of the controlled system. The main aim is to introduce a suitable (scalar)
backward stochastic differential equation (BSDE), with a constraint on the martingale part,
that allows to give a probabilistic representation of the value function of 
the given problem. This is achieved by randomization of control, i.e. by introducing an auxiliary optimization
problem which has the same value as the starting optimal switching problem
and for which the desired BSDE representation is obtained. In comparison with the existing literature
we do not rely on a system of reflected BSDE nor can we use the associated Hamilton-Jacobi-Bellman
equation in our non-Markovian framework.
\end{abstract}

\vspace{5mm}


\vspace{5mm}


\date{}

 \newpage

\section{Introduction}

Stochastic switching control problems arise when a controller acts on a random
system by choosing a piecewise constant control process of the form
$$
\alpha(t)=\xi_0\,1_{[0,\tau_1)}(t)+\sum_{n\ge 1}\xi_n\,1_{[\tau_n,\tau_{n+1})}(t),
\qquad t\ge 0.
$$
Here the switching times $\tau_n$ are an increasing sequence of stopping
times with respect to some given filtration $(\calf_t)$ and the chosen
actions $\xi_n$ are $\calf_{\tau_n}$-measurable random variables
with values in some set $A$, called the set of modes (or regimes).
Thus, when the initial mode $\xi_0=a\in A$ is fixed,
choosing a switching strategy amounts to choosing the double sequence
$
\alpha=(\tau_n,\xi_n)_{n\ge 1}
$. This special form of the strategies is justified when the controller incurs in some cost
whenever the control action is changed, so that only piecewise constant
control processes may have finite cost.
Since optimal
 switching problems are commonly used as models for management issues, they have attracted interest since a long time in the economic
 literature: the interested reader is referred for instance to \cite{Brennanschwarz}, \cite{Dixit} or \cite{DixitPindyck}.

In the classical framework the set of control actions $A$ is finite, say
$A=\{1,\ldots,m\}$. Our main concern is to deal with the case
when the set $A$ is arbitrary which is quite natural for many applications. For instance, each mode $a\in A$ may correspond
to a working regime of a plant, or a production level of a firm;  one may then conceive a situation
when the regime or the production level can be chosen freely within an interval of possible values, 
still retaining the feature that switching from a value to another one entails some cost.\\

Let us now describe our framework. In this paper we will only consider stochastic differential equations in $\R^n$  driven by the Brownian motion.
Suppose initially that the controlled system is described by an equation on the time interval
$[0,T]$ of the form
\beq \label{dynX1_intro}
dX_t^\alpha & = & b( X_t^\alpha, \alpha(t))\,dt +
\sigma( X_t^\alpha, \alpha(t))\,dW_t, \qquad t\in [s,T]\subset [0,T],
\enq
with a given initial condition $X_s^\alpha=x\in\R^n$, where
  $W$ is an  $\R^d$-valued Brownian motion
and the coefficients $b$, $\sigma$ satisfy standard Lipschitz
and growth conditions.
The controller maximizes the reward functional
$$
J(s,x,a,\alpha)   =    \E\Big[\int_s^Tf(X_t^\alpha,\alpha(t))\,dt+g(X^\alpha_T,\alpha(T))
-\sum_{n\ge 1}1_{\tau_n<T}\,c_{\tau_n}(X_{\tau_n}^\alpha,\xi_{n-1},\xi_n)\Big],
$$
where $f$ and $g$ represent the running and terminal rewards
and
$c_t(x,a,a')$ is the cost incurred when switching at time $t$
  from the mode $a$ to the mode $a'$
when the present state  is $x$. The corresponding (so-called) primal value
function of the optimal switching problem, with set of modes equal to $A$, is given at time $s$ by
\begin{equation}\label{primal_ctrolpb_intro}
v_s(x,a)=\sup_{\alpha}J(s,x,a,\alpha).
\end{equation}

Different approaches have been proposed to tackle this problem, that we briefly mention below,
while we refer the reader to \cite{EKb} for a much more detailed discussion.

The classical dynamic programming approach to this problem consists in
studying the associated Hamilton-Jacobi-Bellman equation, which in this case is
a system of partial differential equations coupled by an obstacle condition:
for $a=1,\ldots, m$
\begin{equation}\label{QVIsystem_intro}
\left\{
\begin{array}{l}\dis
\min\left\{
-\partial_sv_s(x,a)-
\call^av_s(x,a) - f(x,a),
v_s(x,a)-\dis\max_{a'\neq a} [v_s(x,a')- c_s(x,a,a')]
\right\}=0,
\\
v_T(x,a)=g(x,a),  \qquad x\in\R^n,\;s\in [0,T],
\end{array}
\right.
\end{equation}
where
$$
\call^av_s(x,a)=\frac12{\,\rm Trace}\,[\sigma(x,a)\sigma(x,a)^TD^2_xv_s(x,a)]
+D_xv_s(x,a)b(x,a)
$$
is the Kolmogorov operator corresponding
to the controlled coefficients $b(x,a)$, $\sigma(x,a)$.

However, such an approach restricts  to the Markovian framework. Among first studies relating the optimal switching problem (with finite number of modes)
 with a system of quasi-variational inequalities of the form (\ref{QVIsystem_intro}) one can cite  \cite{koike}, \cite{LenhartB-83} or 
\cite{TangYong93} and, for general theory concerning stochastic control and its treatment by PDE methods, 
 the interested reader is referred to \cite{BensoussanLions}. 
 More recently, \cite{OlofssonLundstrom} and  \cite{OlofssonLundstrom-2} have further investigated
these systems in the context of filtrations allowing jumps (in that case, the Kolmogorov operator involves an extra non local term). Recent results on
numerical approximation can be found in \cite{Gassiatetal2012} 
for optimal multiple switching problems or in \cite{Bernhard12} for impulse control problems.

Another approach is based on the introduction of a system of Backward
Stochastic Differential Equations (BSDEs). Letting the initial time $s=0$ for simplicity,
one solves a system of reflected BSDEs with interconnected obstacles
looking for unknown adapted processes $(\bar Y^{x,a}_t, \bar Z^{x,a}_t, 
\bar K^{x,a}_t)_{t\in[0,T]}$,
parameterized by $x\in\R^n$ and $a\in A$ and satisfying suitable conditions, such that
\begin{equation}\label{reflected_intro}
\left\{
\begin{array}{l}\dis
\bar Y^{x,a}_t+ \int_t^T \bar Z^{x,a}_s\,dW_s =
g(\bar X^{x,a}_T) + \int_t^Tf(\bar X^{x,a}_s,a)\,ds +\bar K^{x,a}_T-\bar K^{x,a}_t,
\\
\dis \bar Y^{x,a}_t\ge \max_{a'\neq a} [\bar Y^{x,a'}_t- c_t(\bar X^{x,a}_t,a,a')],
\\\dis
\int_0^T\Big[\bar Y^{x,a}_t- \max_{a'\neq a} [\bar Y^{x,a'}_t- c_t(\bar X^{x,a}_t,a,a') ]\Big]\,d\bar K^{x,a}_t=0,
\end{array}
\right.
\end{equation}
where, in particular, $\bar K^{x,a}$ are non decreasing processes, $\bar K^{x,a}_0=0$,
and $\bar X^{x,a}$  are defined by the equations
$$
d\bar X^{x,a}_t =  b( \bar X_t^{x,a}, a)\,dt +
\sigma( \bar X_t^{x,a}, a)\,dW_t, \quad t\in [0,T],
\qquad \; \bar X_0^{x,a}=x.
$$
Under suitable conditions this system is well-posed and one has
a probabilistic representation for the value function: $v_0(x,a)=\bar  Y^{x,a}_0$.

In the two last decades, such a BSDE approach has been extensively used to characterize the primal value function $v_0(x,a)$
(corresponding to (\ref{primal_ctrolpb_intro}) taken at time $s=0$). 
Among the first papers relating the standard optimal switching problem (with $m$ modes) 
to system of reflected BSDEs of the type (\ref{reflected_intro}) one may refer to \cite{Hamjeanblanc}, \cite{HamZh10} or \cite{HuTang10}.
 Some extensions can be found in \cite{Chassaetal}, \cite{HamMor13}, \cite{HamZhao15}, \cite{Fo17}, this list being non exhaustive. 
 In particular, the authors in 
 \cite{HamMor13} and \cite{HamZhao15} combine the BSDE and PDE approach in the Markovian setting where, under appropriate conditions, 
 one can show that  
 the solutions of (\ref{reflected_intro}) and (\ref{QVIsystem_intro}) are related through a standard relation of Feynman-Kac type.


Another approach has been devised, also based on the introduction of a suitable BSDE, but of different type.
Suppose that we are given a Poisson random measure (with finite intensity) on $(0,\infty)\times A$, independent of $W$, and let
  $I^{a}$ denote  the corresponding piecewise constant $A$-valued process starting
from $a\in A$. Let further
$X^{x,a}$ be the solution to
\begin{equation}\label{randomizedsde_intro}
dX^{x,a}_t =  b( X_t^{x,a}, I^{a}_t)\,dt +
\sigma( X_t^{x,a}, I^{a}_t)\,dW_t, \quad t\in [0,T],
\qquad \; X_0^{x,a}=x.
\end{equation}
This will be called the randomized equation, since the switching control
process has been replaced by a random (Poisson) process. Let us then consider
the BSDE
\begin{equation}\label{BSDEconstrained_intro}
\left\{
\begin{array}{l}
 \dis Y_t ^{x,a}
+ \int_t^TZ^{x,a}_s\,dW_s
+\int_{(t,T]}\int_A U_s^{x,a}(a')\,\mu(ds\,da')
  = \ g( X_T^{x,a},  I_T^{a})
  \\\dis\qquad\qquad
+ \int_t^T f(X_s^{x,a},I_s^{a})\, ds + K_T^{x,a} - K_t^{x,a} , \\
\dis U_t^{x,a}(a) \le c_{t-}(X_{t}^{x,a}, I_{t-}^{a},a).
\end{array}\right.
\end{equation}
Here the solution is $(Y_t ^{x,a},Z_t ^{x,a},K_t ^{x,a},U_t^{x,a}(a') )$
($t\in [0,T]$, $a'\in A$)
where the additional martingale term $U$ is a predictable random field
needed to solve the equation with respect to the filtration generated
by the Brownian motion $W$ and the Poisson random measure $\mu$.
This equation is called constrained BSDE, with reference to the inequality required to hold in
\eqref{BSDEconstrained_intro}.
Under suitable assumptions there exists a unique minimal solution
(in a sense to be defined) and
it is proved that the value function is also represented by the
formula $v_0(x,a)=Y_0^{x,a}$.
This  control randomization method was introduced in
\cite{bou09}. There the author also formulates a corresponding randomized
optimal control problem (i.e. an auxiliary or dual problem)
and a stochastic target problem related to optimal switching.
In the framework of switching problems and associated BSDEs
the method was further developed and extended in
\cite{EKa}, \cite{EKb}, \cite{EKc} 
and later applied to different contexts by many authors,
see for instance  
 \cite{KMPZ10},  \cite{KP12},
 \cite{BCFP16bAAP},  \cite{BCFP16bSPA}
  \cite{FP15},  \cite{CossoFuhrmanPham16}, \cite{CoCoFu18},
\cite{Fo17}, \cite{FuPhZe16}, \cite{Bandini15},
\cite{BandiniFuhrman15}.

We note that the two approaches based on BSDEs have immediate generalization,
which already appear in many of the references cited above,
to the case of path-dependent coefficients (also called the non-Markovian case),
that is when the value at time $t$ of the drift and the diffusion also depend on the past
history $(X_s^\alpha)_{s\in [0,t]}$ of the controlled process.
Moreover the approach based on the contrained BSDE is more promising from a computational point of view since one
has to deal with a single equation instead of a system: as such, numerical methods have been devised to treat this equation,
see \cite{KLP14}, \cite{KLP15}.

Finally, we cite another special approach to optimal switching developed in \cite{DjHamPo09}, 
which works
both in Markovian and non-Markovian situations, 
where BSDEs are replaced by an implicit optimal stopping problem.

As mentioned above, our main concern in this paper is to address the switching problem when the set $A$ is infinite (not necessarily countable). For greater generality 
we will consider path-dependent coefficients and try to generalize the approaches based on BSDEs. While addressing an infinite system of reflected BSDEs of the form \eqref{reflected_intro}, or using the approach 
of \cite{DjHamPo09}, seems difficult, it turns out that a generalization of the approach based on the constrained 
BSDE is possible, and this is in fact the main content of the present paper.
Another motivation is the fact that we will still be concerned with a single BSDE
even if the number of modes is infinite, so the feasibility
of numerical approximation will be preserved, although we will not deal
with this issue 
in this paper.

Following \cite{bou09}  and \cite{EKb}, we introduce 
an auxiliary optimization problem, called randomized control problem
(see section \ref{Randomized}
for a precise formulation),
having the same value as the original
switching problem  and we show that this common value can be represented by means of the solution to the constrained
BSDE \eqref{BSDEconstrained_intro}, even when the set of modes $A$ is infinite. To this aim
we have to find entirely new proofs. Indeed, in \cite{bou09} 
the result was proved by showing that the switching problem and
the randomized problem correspond to the same Hamilton-Jacobi-Bellman
equation, since in that paper only the Markovian case was addressed.
On the contrary in \cite{EKb} the non-Markovian situation was studied, but
the link between the randomized problem and the switching problem was 
proved by means of the system of reflected BSDEs \eqref{reflected_intro}, which does not seem easy to solve in the case when $A$ is infinite. 
In fact, we establish a direct link between the switching
problem and the randomized one, and between the latter and the constrained BSDE 
\eqref{BSDEconstrained_intro}. As a consequence our treatment is almost entirely self-contained,
except for some technical results related to the randomization technique.

A drawback of the randomization method is that it does not immediately
provide a description nor even the existence of an optimal control, but it rather aims
at a convenient representation of the value function. However, it also works in situation where
an optimal control may not exist, for instance without compactness
assumptions on the set of modes $A$.  

The model that we formulate for the switching problem is fairly general:
all coefficients, including the switching costs, are path-dependent
and may be unbounded. 
On the diffusion coefficient (the volatility), that can also be controlled, 
we do not impose any nondegeneracy condition which implies that the case of
deterministic
optimal switching falls under the scope of our results.

To complete our discussion on the possible approaches to optimal
switching with infinitely many modes
we finally mention that results based on Hamilton-Jacobi-Bellman
equations have been obtained, but limited to the Markovian case when
there is no path-dependence in the coefficients.
In fact in this case optimal switching can be considered as a special case
of optimal impulse problem, where the state of the system is the pair $(X_t,I_t)$.
The randomization method has been successfully used in this context as well
in \cite{KMPZ10}. However from a technical point of view these results are not always satisfactory
since they impose stringent assumptions, being designed to hold true
in a more general (or simply different) context. We believe that building on our approach
more refined results can be obtained in the case of Markovian optimal switching with infinitely many modes, and
this will be the object of future research.
On the contrary, there are not many results on optimal impulse control
in the non-Markovian context that apply to general
models; one example is \cite{DjHamHdh10}, which however seems difficult
to generalize to the case of an infinite set $A$.

The plan of the paper is as follows: in section 
\ref{S:Formulation} we formulate our assumptions and introduce
the optimal switching problem.
In section 
\ref{Randomized} we formulate the auxiliary  randomized optimization problem 
and prove that its value coincides with the value of the optimal
switching problem. The proof is rather technical and is presented in
section
\ref{proofthm}. Finally in section
\ref{Sec:separandom} we show that a constrained BSDE of the form
\eqref{BSDEconstrained_intro}
can be associated to the randomized problem thus giving
the desired representation of the value for the starting optimal switching
problem as well.

\section{General assumptions and formulation of the optimal switching problem}
\label{S:Formulation}

\subsection{General notations  and assumptions}
\label{SubS:Notation}

We start this section by an informal description of our optimization problem.

In the following we will consider controlled stochastic equations in $\R^n$ of the form
\beq \label{dynX1}
dX_t^\alpha & = & b_t( X^\alpha, \alpha(t))\,dt +
\sigma_t( X^\alpha, \alpha(t))\,dW_t,
\enq
for $t\in [0,T]$, where $T>0$ is a fixed deterministic and finite
terminal time, and an initial condition  $X_0^\alpha=x_0$, a given deterministic
point in $\R^n$. $W$ is a standard Brownian motion with values in $\R^d$.
The control process $\alpha(\cdot)$ is a switching process:
it takes values in a set $A$, called set of modes (or regimes),
and it is piecewise constant: it starts at a deterministic
mode $\xi_0\in A$ and at random jumps times $\tau_n$
it jumps from  $\xi_{n-1}$ to $\xi_n$ ($n\ge 1$).
$\tau_n$ are stopping times for the filtration $(\calf^W_t)$ generated by $W$
and modes $\xi_n$ are also random $A$-valued variables,
each assumed to be $\calf^W_{\tau_n}$-measurable.

In our framework we include
 path-dependent (or hereditary) systems, i.e. exhibiting memory
effects with respect to the state. Indeed, the value of the coefficients
$b,\sigma$  at time $t$ depend
on the values $X_s^\alpha$ for $s\in [0,t]$: this non-anticipative dependence
will be expressed below in a standard way
by requiring that the coefficients should be progressive with respect to
the canonical filtration on the space of continuous paths.

The reward functional, to be maximized over an appropriate class
of switching processes $\alpha$, has the form
$J(\alpha) $ $ =$ $   J_1(\alpha)-J_2(\alpha)$,
where
$$
J_1(\alpha)   =    \E\Big[\int_0^Tf_t(X^\alpha,\alpha(t))\,dt+g(X^\alpha,\alpha(T))\Big],
\qquad
J_2(\alpha)    =    \E\Big[\sum_{n\ge 1}1_{\tau_n<T}\,c_{\tau_n}(X^\alpha,\xi_{n-1},\xi_n)\Big].
$$
The functional $J_1$ has a classical form, and also contains real-valued path-dependent
coefficients $f,g$;  the functional $J_2$ takes into account the cost of switching:
the (path-dependent) nonnegative function
$c_t(x,a,a')$ is interpreted as the cost incurred when switching at time $t$
  from the mode $a$ to the mode $a'$
when the  trajectory is $x(\cdot)$.

\bigskip

Now let us introduce notations and precise assumptions  on the data
$A,b,\sigma,f,g,c,x_0,\xi_0$.
In the next paragraph we will formulate   the optimization problem
by describing in particular the class of admissible switching strategies.

Let us denote by $\bfC_n$   the  space  of continuous paths from $[0,T]$ to $\R^n$, equipped  with the usual supremum norm
$\|x\|_{_\infty}$ $=$ $x^*_T$, where we set  $x^*_t$ $:=$ $\sup_{s\in [0,t]}|x(s)|$, for $t$ $\in$ $[0,T]$ and $x$ $\in$ $\bfC_n $. We introduce the filtration
$(\calc_t^n)_{t\in[ 0,T]}$, where by $\calc_t^n$ we denote the $\sigma$-algebra generated by the canonical coordinate maps $\bfC_n \to \R^n$,
$x(\cdot)\mapsto x(s)$ up to time $t$, namely
\beqs
\calc_t^n  &:=  &  \sigma \{ x(\cdot)\mapsto x(s)\; :\, s\in [0,t] \},
\enqs
and we denote $Prog(\bfC_n )$ and $\calp(\bfC_n )$  the progressive and
predictable
 $\sigma$-algebra on $[0,T]\times\bfC_n $ with respect to $(\calc_t^n)$,
 respectively. [Indeed, one can prove that these  $\sigma$-algebras essentially coincide:
 see Remark (8.4) in Chapter V of \cite{rowi},
 but we will not need this for the sequel.]

We require 
the space of control actions $A$ to be a  Borel space. We  recall that a  Borel  space    is a topological space  homeomorphic to a  Borel subset of a Polish space.
The terminology Lusin space, instead of Borel space, is sometimes used.
The space $A$ will be  endowed with its Borel $\sigma$-algebra  $\calb(A)$.

Throughout the paper, the following assumptions will be in force.

\vspace{3mm}

\noindent {\bf (A1)}
\begin{itemize}
\item [(i)]  $A$ is a  Borel space.
\item [(ii)]
The functions $b,\sigma,f$ are defined on $[0,T]\times \bfC_n\times   A$ with
values in $\R^n$, $\R^{n\times d }$ and $\R$ respectively,  they are assumed
to be $Prog(\bfC_n) \otimes \calb(A)$-measurable
(see also Remark \ref{A1_misurabilita} below).
\newline
The function $c$ is defined on
$[0,T]\times \bfC_n\times   A\times A$, it takes nonnegative real values and it is
assumed to be $\calp(\bfC_n) \otimes \calb(A)\otimes \calb(A)$-measurable.
\newline
The function $g$ is defined on
$\bfC_n\times  A$ and takes real values.

\item [(iii)] For every $t\in[0,T]$, the functions
  $g(x,a)$, $ b_t(x,a)$,  $\sigma_t(x,a)$ and $f_t(x,a)$
are
 continuous functions of  $(x,a)\in\bfC_n\times A$
 ($\bfC_n$ being equipped with the supremum norm).
 \newline
The
 function
  $c_t(x,a,a')$
is a
 continuous functions of  $(t,x,a,a')\in[0,T]\times\bfC_n\times A\times A$.

\item [(iv)]
 There exist nonnegative constants  $L$ and $r$ such that
\beq
|b_t(x,a) - b_t(x',a)| + |\sigma_t(x,a)-\sigma_t(x',a)|
& \leq & L   (x-x')^*_t, \label{lipbsig} \\
|b_t(0,a)| + |\sigma_t(0,a)| & \leq & L,  \label{borbsig}
\\
|f_t(x,a)| + |g(x,a)| + |c_t(x,a,a')| & \leq & L \big(1 + (x^*_t)^r \big), \label{PolGrowth_f_g}
\enq
for all $(t,x,x',a,a')$ $\in$ $[0,T]\times\bfC_n \times\bfC_n \times A\times A$.
\item [(v)]  $x_0\in\R^n$ and $\xi_0\in A$ are given: they represent the initial
state and mode, respectively.
\end{itemize}

\vspace{2mm}

\begin{Remark} \label{A1_misurabilita}
{\rm
The measurability conditions in
{\bf (A1)}-(ii) entail   the following property,
which is easily verified:
\begin{itemize}

\item [(ii)'] Whenever
$(\Omega,\calf,\P)$ is a probability space with a filtration $\F$,  $X$
 is an $\F$-progressive process
 with values in $\R^n$, and $a,a'\in A$, then the processes
 $b_t( X , a)$,
$\sigma_t( X , a)$,
$f_t(X,a)$, $c_t(X,a,a')$, defined for ${t\in [0,T]}$, are also
   $\F$-progressive.
\end{itemize}
All the results in this paper still hold, without any change in the proofs,
if property (ii)' is assumed to hold instead of (ii). In some cases, (ii)' is
easier to be checked directly.

We finally note that the function $g$, being continuous, is also Borel measurable
(equivalently, it is $\calc_T^n  \otimes \calb(A)$-measurable).
\ep
}
\end{Remark}

\begin{Remark}\label{degeneracy}
{\rm  We mention that no non-degeneracy assumption on the diffusion coefficient $\sigma$ is imposed.
 In particular the case of deterministic switching, where $\sigma=0$, is included, and in this special
 case there is of course no need to introduce a Wiener process nor a probability space.
\ep
}
\end{Remark}


\subsection{Formulation of the optimal switching  problem}
\label{Primal}

We assume that
$A,b,\sigma,f,g,c,x_0,\xi_0$
 are given and
satisfy the assumptions {\bf (A1)}. A setting
$(\Omega$, $\calf$, $\P$, $W)$
for the optimization
problem consists of  a complete probability space $(\Omega,\calf,\P)$ and
an $\R^d$-valued  process
 $W$ which is a standard Wiener process with respect to $\P$.

Let us denote $\F^W=(\calf^W_t)_{t\ge 0}$ the right-continuous and $\P$-complete filtration generated by $W$. We define
the set $\Ac$ of admissible control strategies: its elements are the
double sequences of the form
$$
\alpha=(\tau_n,\xi_n)_{n\ge 1},
$$
where:
\begin{enumerate}
\item[(i)] each $\tau_n$  is an $\F^W$-stopping time;
\item[(ii)] each $\tau_n$   takes values in $(0,\infty]$ and the sequence
$(\tau_n)_{n\ge 1}$ is nondecreasing, a.s.;

\item[(iii)]  if $\tau_n<\infty$ then $\tau_n<\tau_{n+1}$, for every $n\ge 1$, a.s.;
\item[(iv)] each $\xi_n$ is a random variable
with values in $A$, which is $\calf^W_{\tau_n}$-measurable;
\item[(v)] $\tau_n\to \infty$ a.s.
and $\tau_n\neq T$ a.s. for every $n\ge 1$.
\end{enumerate}

\begin{Remark}\label{strategiecomemisure}
{\rm  Conditions $(i)-(iv)$ can be restated by saying that
$\alpha$ is a marked (or multivariate)
point process in $A$. It is convenient in the following
to use this definition although the control horizon
$T$ is finite. The condition $\tau_n\to\infty$
 can be expressed by saying that the explosion time
$\lim_n\tau_n$ is infinite a.s.
We comment further on this condition and on the condition $\tau_n\neq T$
in Remark
\ref{varianti e condizioni}.
\ep
}
\end{Remark}

Given $\alpha\in\cala$, we introduce the associated piecewise constant
process, denoted by $\alpha(\cdot)$ (with a slight abuse of notation)
and defined as
$$
\alpha(t)=\xi_0\,1_{[0,\tau_1)}(t)+\sum_{n\ge 1}\xi_n\,1_{[\tau_n,\tau_{n+1})}(t),
\qquad t\in [0,T],
$$
where $\xi_0$ is the given starting mode. Notice that the formal sum makes
obvious sense
even if there is no addition operation defined in $A$.

The corresponding trajectory $X^\alpha$ is defined as the solution
to the controlled equation
\begin{equation}\label{stateeq}
    dX_t^\alpha \ = \ b_t( X^\alpha, \alpha(t))\,dt +
\sigma_t( X^\alpha, \alpha(t))\,dW_t
\end{equation}
on the interval $[0,T]$ with initial condition $X_0^\alpha=x_0$.
Since we assume that {\bf (A1)} holds,  by standard results (see e.g. \cite{rowi} Thm V. 11.2,
 or   \cite{jacod_book}  Theorem 14.23),
there exists an almost surely
unique $\F$-adapted strong solution $X^\alpha$ $=$ $(X_t^\alpha)_{t\in[0, T]}$   to \eqref{stateeq}
with continuous trajectories a.s. and such that
\begin{equation}\label{growthsolprimal}
        \E\,\Big[\sup_{t\in [0,T]}|X_t^\alpha|^p\Big]  \le C_p <  \infty,
\end{equation}
for every $p\in [1,\infty)$, where the constant $C_p$, depends only on $p,T,n,d$ and
the constants $L,r$ appearing in Assumption
{\bf (A1)}.
The stochastic optimal control problem  under partial observation consists in maximizing, over all $\alpha\in\Ac$, the reward functional
\begin{equation}\label{gaineq}
J(\alpha)  =   J_1(\alpha)-J_2(\alpha),
\end{equation}
where
\beq \label{Juno}
J_1(\alpha) &  = &  \E\Big[\int_0^Tf_t(X^\alpha,\alpha(t))\,dt+g(X^\alpha,\alpha(T))\Big],
\\\label{Jdue}
J_2(\alpha) &  = &  \E\Big[\sum_{n\ge 1}1_{\tau_n<T}\,c_{\tau_n}(X^\alpha,\xi_{n-1},\xi_n)\Big].
\enq
We define the value of the optimal switching problem as
\begin{equation}\label{primalvalue}
{\text{\Large$\upsilon$}}_0 = \sup_{\alpha\in\Ac} J(\alpha).
\end{equation}

Since we do not impose growth conditions on the cost
function $c$, it is possible that $J_2(\alpha)=\infty$ for some admissible
$\alpha\in\cala$. However, we have the following simple result.

\begin{Lemma}
There exists a finite constant $C$, depending only on $T,n,d$ and
the constants $L,r$ appearing in assumptions
{\bf (A1)}, such that $|{\text{\Large$\upsilon$}}_0 |\le C$.
\end{Lemma}

 \noindent {\bf Proof.} By standard estimates on the state equation
 (the same ones leading to \eqref{growthsolprimal}) and the growth conditions
 imposed in  {\bf (A1)}-(iv) it is easily shown that $|J_1(\alpha)|\le C$
 for every $\alpha\in\cala$.
 Since $c$ is nonnegative we have $J(\alpha)\le J_1(\alpha)\le C$
 for $\alpha\in\cala$ and it follows that  ${\text{\Large$\upsilon$}}_0\le C$.

 Now let us consider the strategy $\bar\alpha$ without switchings
 (i.e. such that $\tau_n=\infty$ for $n\ge 1$). Then we have $J_2(\bar\alpha)=0$
 and so
$$ {\text{\Large$\upsilon$}}_0\ge J(\bar\alpha)=J_1(\bar\alpha)\ge -C,
$$
and we conclude that $|{\text{\Large$\upsilon$}}_0 |\le C$.
\qed

We end this section with several comments on the previous formulation
of the optimization problem and its possible variants.

\begin{Remark} \label{varianti e condizioni}
{\rm
\begin{enumerate}
\item According to large part of the literature on optimal
switching, we do not allow for a switching at initial time $t=0$.
This is not a real loss of generality, since a switching at time $0$
does not affect the controlled trajectory  $X^\alpha$ and
it is easy to reduce the problem to the formulation that we adopt.
\item  In our definition of admissible strategy we have imposed
the condition of being non-explosive. This implies
 that
 $$N_T:=\sum_{n\ge 1}1_{\tau_n\le T}
 $$
 is finite a.s., meaning that infinitely many
switchings in the time interval within the control horizon $T$
are not allowed. Alternatively, one
may impose that there exists $\delta>0$ such that
$c_t(x,a,a')\ge \delta$ for every $t\in [0,T]$, $x\in \bfC_n$,
$a,a'\in A$, which is a common requirement in the literature
on switching problems. Under this additional assumption,
any strategy $\alpha$ with $N_T=\infty$ has $J(\alpha)=-\infty$
and cannot be optimal. We will not need that
$c_t(x,a,a')\ge \delta$ and will only assume the weaker
conditions that $c_t(x,a,a')\ge 0$ and $\tau_n\to\infty$ for every admissible strategy.

\item Often, the following assumption is imposed on the cost function:
for every distinct
$a_1,a_2,a_3\in A$ and
for every $t\in [0,T]$, $x\in \bfC_n$,
\begin{equation}\label{triangle}
    c_t(x,a_1,a_3)<c_t(x,a_1,a_2)+c_t(x,a_2,a_3).
\end{equation}
This says that switching from mode $a_1$ to mode $a_3$ directly
is more convenient than a double switching from mode
 $a_1$ to $a_2$ followed immediately by a  switching from $a_2$ to  $a_3$.
 This condition entails that any strategy
 for which a switching time $\tau_n$ equals $\tau_{n+1}$ cannot
 be optimal. We will not need the condition \eqref{triangle},
 but we have imposed that $\tau_n<\tau_{n+1}$ (whenever
 $\tau_n$  is finite).

 \item A variant of the optimal switching problem is obtained by
 allowing for a switching at the terminal time, that is by removing
 the requirement that $\tau_n\neq T$ and
 modifying the functional $J_2$, introduced in \eqref{Jdue}, in the following way:
\begin{equation}\label{J2modificato}
J_2(\alpha)    =
\E\Big[\sum_{n\ge 1}1_{\tau_n\le T}\,c_{\tau_n}(X^\alpha,\xi_{n-1},\xi_n)\Big],
\end{equation}
in order to
take into account the cost of a switching at the final time.
In some papers, the following condition is imposed on the data:
for every
$a\in A$ and
 $x\in \bfC_n$,
\begin{equation}\label{finalswitch}
    g(x,a)> \sup_{a'\in A,a'\neq a}(g(x,a')-c_T(x,a,a')).
\end{equation}
This says that at the final time it is more convenient to remain in the current
mode $a$ rather than switching to any another mode $a'$, which would
give a reward $g(x,a')$ but would incur in a cost $c_T(x,a,a')$.  If
\eqref{finalswitch} is required, the optimization problem has the same
value (and the same optimal control, if it exists) whether $J_2$
is defined by
\eqref{Jdue} or by
\eqref{J2modificato}.

In this paper we will not impose condition \eqref{finalswitch}
but we require that $\tau_n\neq T$ a.s.

\item In some papers a different, weak formulation
of the optimization problem is considered. Sometimes
this leads to additional assumptions (for instance in
 \cite{HuTang10}  $\sigma$ does not depend
on the control and it is assumed to be invertible) or to a more
involved formulation where the probability space
and the Wiener process are also part of the control
(as in \cite{FP15}). In this paper we have chosen 
the simpler and more natural strong formulation.
We mention that whenever ultimately  the value of the
optimization problem  
is represented by means of a uniquely solvable
BSDE then obviously the value remains the same
for both formulations.

\end{enumerate}
\ep
}
\end{Remark}

\section{The randomized stochastic optimal control problem}
\label{Randomized}

We still assume that
$A,b,\sigma,f,g,c,x_0,\xi_0$
 are given and
satisfy the assumptions {\bf (A1)}.
We introduce an auxiliary optimization problem,
 that we call  randomized  optimal control  problem, and we will eventually prove
 that it has the same value as the optimal switching problem
 formulated in section \ref{Primal}.
However, the randomized problem has the advantage that it can be directly
related to a suitable  BSDE, as we will see in the following sections.

 To this end we need one additional datum, that
will play the role of an intensity measure for a Poisson process:

\vspace{3mm}

\noindent {\bf (A2)}
Let  $\lambda$ be a finite positive measure on  $(A,\calb(A))$
with full topological support.

\vspace{3mm}

Since $A$ is separable (as a Borel space), such a measure
always exists: for instance, one could choose a convex
  combination of Dirac measures at points $a_i\in A$,
  where $(a_i)$ is a dense sequence in $A$. In general
there are many possible choices for the measure $\lambda$
 and in any case {\bf (A2)} is not a restriction imposed on the
original optimization problem. It will be assumed to hold from now on.

\subsection {Formulation of the randomized control problem}
\label{randomizedformulation}

We say that $(\hat \Omega, \hat \calf,\hat \P,  \hat W, \hat \mu)$
is a setting for the randomized control problem if
$(\hat \Omega, \hat \calf,\hat \P)$
is an arbitrary complete probability space,
the  process
$\hat W$
is a standard Wiener process in $\R^{ d}$ under $\hat \P$,
$\hat \mu$ is a Poisson random measure on  $(0, \infty) \times A$ with intensity $\lambda(da)dt$
  under $\hat \P$, independent of $W$.
Thus, $\hat \mu$ is a sum
 of random Dirac measures and it has
 the form $\hat\mu=\sum_{n\ge 1}\delta_{(\hat \sigma_n,\hat \eta_n)}$,
 where $(\hat \eta_n)_{n\ge 1}$ is a sequence of $A$-valued random variables
 and
 $(\hat \sigma_n)_{n\ge 1}$ is a strictly increasing sequence of random variables
with values in $(0,\infty)$, and for any $C\in\calb(A)$ the process
$\hat \mu((0,t]\times C)-t\lambda(C)$, $t\ge 0$, is a $\hat \P$-martingale.
We also define the piecewise-constant $A$-valued process associated to $\hat{\mu}$
and starting at the initial mode $\xi_0$:
\begin{equation}
\label{I}
\hat I_t \ =
\xi_0\,1_{[0,\hat \sigma_{ 1})}(t)+
\sum_{n\ge 1}\hat \eta_n\,1_{[\hat \sigma_n,\hat \sigma_{n+1})}(t), \qquad t\ge 0.
\end{equation}
The formal sum in \eqref{I} makes sense
even if there is no addition operation defined in $A$, but
when $A$ is a subset of a linear space  formula \eqref{I} can be written as
\[
\hat I_t \ = \xi_0 + \int_0^t\int_A(a-\hat I_{s-})\,\hat \mu(ds\, da), \qquad t\ge 0.
\]
Let $\hat X$ be the solution to the equation
\beq \label{dynXrandom}
d\hat X_t &=&  b_t( \hat X,\hat I_t)\,dt + \sigma_t(\hat X,\hat I_t)\,d\hat{W}_t,
\enq
for $t\in [0,T]$, starting from $\hat X_0$ $=$ $ x_0$, the initial state fixed
at the beginning.

We introduce the filtration
$\F^{\hat W,\hat \mu}=(\calf^{\hat W,\hat \mu}_t)_{t\ge 0}$
generated by  $\hat W,\hat \mu$ and defined  by the formula:
\beq\label{filtraznaturaleWmu}
\calf^{\hat W,\hat \mu}_t&=&\sigma (\hat W_s,
\hat \mu((0,s]\times C)\,:\, s\in [0,t],\, C\in\calb(A))
\vee \caln,
\enq
where $\caln$ denotes the family of $\hat \P$-null sets of $\hat \calf$.
We denote $\calp(\F^{\hat W,\hat \mu})$ the corresponding predictable
$\sigma$-algebra.

Under {\bf (A1)} it is well-known (see e.g. Theorem 14.23 in \cite{jacod_book}) that there exists
an almost surely unique $\F^{\hat W,\hat \mu}$-adapted strong solution
$\hat X$ $=$ $(\hat X_t)_{t\in [0, T]}$   to \eqref{dynXrandom}, satisfying
$\hat X_0= x_0$, with continuous trajectories a.s. and such that for every $p\in[1,\infty)$,
\beq\label{EstimateX}
    \hat \E\,\Big[\sup_{t\in [0,T]}|\hat X_t|^p\Big]  & \le & C_p,
\enq
where $C_p$ is a finite constant whose value depends only on $p$, $T$, $n$, $d$
and the constants $L,r$ occurring in {\bf (A1)}-(iv).

We can now formulate the randomized optimal control problem as follows.
We introduce the set $\hat \calv$
of admissible controls as the set of all $\hat\nu=\hat \nu_t(\hat \omega,a):
\hat \Omega\times \R_+\times A\to (0,\infty)$,
which are $\calp(\F^{\hat W,\hat \mu})\otimes \calb(A)$-measurable and bounded.
To any $\hat \nu$ in $\hat \calv$, we associate its Dol\'eans-Dade exponential process $\kappa_t^{\hat\nu} $ defined as follows
\beq \nonumber
\kappa_t^{\hat\nu} \ &=& \
\Ec_t\bigg(\int_0^\cdot\int_A (\hat \nu_s(a) - 1)\,(\hat\mu(ds\, da)- \lambda(da)\,ds)\bigg) \\
&=& \ \exp\left(\int_0^t\int_A (1 - \hat \nu_s(a))\lambda(da)\,ds
\right)\prod_{0<\hat \sigma_n\le t}\hat \nu_{\hat \sigma_n}(\hat \eta_n),\qquad t\ge 0.\label{doleans}
\enq
It is known that $\kappa^{\hat \nu}$ is a martingale with respect to $\hat\P$ and $\F^{\hat W,\hat \mu}$
and thus we define a new probability measure by setting
$\hat\P^{\hat\nu}(d\hat\omega)=\kappa_T^{\hat\nu}(\hat\omega)\,\hat\P(d\hat\omega)$. From the
Girsanov theorem for multivariate point processes (\cite{ja}) it follows that under $\hat \P^{\hat\nu}$
the $\F^{\hat W,\hat \mu}$-compensator of $\hat\mu$ on the set
$[0,T]\times A$ is the random measure $\hat\nu_t(a)\lambda(da)dt$. Moreover, $\hat W$ remains a
standard Wiener process under $\hat\P^{\hat\nu}$, so that using both Assumptions \eqref{lipbsig}-\eqref{borbsig} and standard results
we obtain the following generalization of the estimate \eqref{EstimateX}:
\beq\label{EstimateX_nu}
    \sup_{\hat\nu\in \hat \Vc}\,\hat\E^{\hat\nu}\,\Big[\sup_{t\in [0,T]}|\hat X_t|^p\Big]  & \le & C_p,
\enq
where $\hat\E^{\hat\nu}$ denotes the expectation with respect to $\hat\P^{\hat\nu}$
and $C_p$ is the same as in \eqref{EstimateX}.
We finally introduce the reward functional of the randomized control problem

\beq \label{defJrandomized}
J^\Rc(\hat\nu) &=&  J_1^\Rc(\hat\nu)-J_2^\Rc(\hat\nu),
\enq
where
\beq \label{JRuno}
J_1^\Rc(\hat\nu) &=&  \hat\E^{\hat\nu}
\Big[\int_0^Tf_t(\hat X,\hat I_t)\,dt+g(\hat X, \hat I_T)\Big],
\\\label{JRdue}
J_2^\Rc(\hat\nu) &  =
&  \hat\E^{\hat\nu}
\Big[\sum_{n\ge 1}1_{\hat\sigma_n<T}
\,c_{\hat \sigma_n}(\hat X,\hat \eta_{n-1},\hat\eta_n)\Big],
\enq
where we use the convention $\hat \eta_0=\xi_0$.  We note that
 $$
J^\Rc(\hat\nu ) =  \hat\E
\left[\kappa_T^{\hat\nu }\left( \int_0^Tf_t(\hat X,\hat I_t)\,dt+g(\hat X,\hat I_T)
-
\sum_{n\ge 1}1_{\hat\sigma_n<T}\,c_{\hat\sigma_n}(\hat X,\hat\eta_{n-1},\hat\eta_n)
\right)\right]
 $$
is always finite: indeed, letting $\hat N_T=\sum_{n\ge 1}1_{\hat\sigma_n\le T}$
and recalling the growth conditions
in  \eqref{PolGrowth_f_g} we see that
\begin{equation}\label{stimasukappa}
0\le \kappa_T^{\hat\nu }\le 
\exp(T\lambda(A)\,(1+\sup \hat \nu))\cdot (\sup\hat \nu)^{\hat N_T},
\end{equation}
\begin{equation}\label{stimasureward}
\left| \int_0^Tf_t(\hat X,\hat I_t)\,dt+g(\hat X,\hat I_T)\right|
+
\sum_{n\ge 1}1_{\hat\sigma_n<T}\,c_{\hat\sigma_n}(\hat X,\hat\eta_{n-1},\hat\eta_n)
\le C (1+\sup_{t\in [0,T]}|\hat X_t|)^r \, (1+\hat N_T),
\end{equation}
for a suitable constant $C$. Noting that $\hat N_T$ has Poisson law
with parameter $\lambda(A)\,T$ under $\hat \P$ and recalling \eqref{EstimateX},
we see that the right-hand sides in the above expressions lie in
$L^p(\hat\P)$ for every $p\in [1,\infty)$ and the finiteness
of
 $J^\Rc(\hat\nu)$ follows.

The randomized stochastic optimal control problem consists in maximizing
$ J^\Rc(\hat\nu)$ over all $\hat\nu\in\hat\calv$. Its value
is defined as
\begin{equation}\label{dualvalue}
{\text{\Large$\upsilon$}}_0^\Rc \;=\;   \sup_{\hat \nu\in\hat \calv} J^\Rc(\hat \nu).
\end{equation}

\begin{Remark}\label{infzero}\emph{A comparison between the 
starting optimal switching problem and the randomized problem
may be useful. In the switching problem, the switching process $\alpha(\cdot)$ is chosen to control the system.
In the randomized problem $\alpha(\cdot)$ is first 
replaced by the Poisson point process $\hat I(\cdot)$ (associated with random measure $\hat \mu$) in the coefficients
of the equation solved by $\hat X$.
In this new problem, the effect of a control strategy $\hat \nu $ is to modify the intensity of $\hat I$ 
(more precisely, to change its compensator from $\lambda(da)dt$
to $\hat\nu_t(a)\lambda(da)dt$) and thus  also to affect the 
law of the process $\hat X$. 
This is done by introducing the probabilities
$\hat\P^{\hat \nu}$ via the Girsanov theorem, and optimizing the reward functional $\mathcal{J}^{\cal R}(\hat \nu)$
among this family of equivalent probability measures parameterized by the set of all bounded predictable random fields $\hat\nu$.
}
\ep
\end{Remark}

\begin{Remark}\label{infzero}\emph{
Let us define
$\hat\Vc_{\inf\,>\,0}=\{\hat\nu\in \hat\Vc\,:\, \inf_{\hat \Omega\times[0,T]\times A}\hat\nu>0\}$.
Then we claim that
\begin{equation}\label{eqinfzero}
{\text{\Large$\upsilon$}}_0^\Rc \;= \;
\sup_{\hat\nu\in\hat\Vc_{\inf\,>\,0}}J^\Rc(\hat\nu).
\end{equation}
Indeed, given $\hat\nu\in\hat\Vc$ and $\epsilon>0$, define $\hat\nu^\epsilon=
\hat\nu\vee \epsilon\in \hat\Vc_{\inf\,>\,0}$
and write the gain \eqref{defJrandomized} in the form
 $$
J^\Rc(\hat\nu^\epsilon) =  \hat\E
\left[\kappa_T^{\hat\nu^\epsilon}\left( \int_0^Tf_t(\hat X,\hat I_t)\,dt+g(\hat X,\hat I_T)
-
\sum_{n\ge 1}1_{\hat\sigma_n<T}\,c_{\hat\sigma_n}(\hat X,\hat\eta_{n-1},\hat\eta_n)
\right)\right].
 $$
 As noted earlier, the expression in curve brackets lies in
 $L^p(\hat\P)$ for every $p\in [1,\infty)$. Moreover we have
$ \kappa_T^{\hat\nu^\epsilon}\to \kappa_T^{\hat\nu}$ a.s. as $\epsilon\to 0$,
and using the estimate \eqref{stimasukappa} with $\nu^\epsilon$ instead of
$\nu$ we conclude that $ \kappa_T^{\hat\nu^\epsilon}\to \kappa_T^{\hat\nu}$ in
 $L^p(\hat\P)$ for every $p\in [1,\infty)$ as well.
 It follows that  $J^\Rc(\hat\nu^\epsilon) \to J^\Rc(\hat\nu) $,
 which implies
 $${\text{\Large$\upsilon$}}_0^\Rc=\sup_{\hat\nu\in\hat\Vc}J^\Rc(\hat\nu)\le \sup_{\hat\nu\in\hat\Vc_{\inf\,>\,0}}J^\Rc(\hat\nu).
 $$
 The other inequality being obvious, we obtain  \eqref{eqinfzero}.
}
\ep
\end{Remark}

\begin{Remark}\label{indepofthesetting}\emph{
We stress the fact that the value ${\text{\Large$\upsilon$}}_0^\Rc$
of the randomized control problem
defined in \eqref{dualvalue}
does not depend on the specific
setting $(\hat \Omega, \hat \calf,\hat \P,  \hat W, \hat \mu)$ that is chosen
in its formulation.
\newline
More precisely, this means that if
 $(\tilde \Omega, \tilde \calf,\tilde \P,  \tilde W, \tilde \mu)$ is another setting with
 the properties described at the beginning of this section, and
 if  the corresponding value $\tilde{\text{\Large$\upsilon$}}_0^\Rc$ is  defined in analogy
 with what was done before then we have the equality
 ${\text{\Large$\upsilon$}}_0^\Rc$ $=$ $\tilde{\text{\Large$\upsilon$}}_0^\Rc$.
\newline
We do not write down the proof of this statement, since it is entirely analogous to
 Proposition 3.1 of \cite{BCFP16bAAP}, where a classical optimization problem
 with continuous control was addressed instead of a switching problem,
 but the arguments remain the same.
\newline
As a consequence, we obtain the rather intuitive conclusion that
the value
${\text{\Large$\upsilon$}}_0^\Rc$
  is just a functional of the (deterministic) elements
$A,b,\sigma,f,g,c,x_0,\xi_0,\lambda$ appearing
in the assumptions {\bf (A1)} and {\bf (A2)}.
 Later on, in
 Theorem
\ref{MainThm},  we will prove
that in fact ${\text{\Large$\upsilon$}}_0^\Rc$
does not depend on the choice of $\lambda$  either.
}
\ep
\end{Remark}

\begin{Remark}\label{productextension}
\emph{
Starting
from a setting  $(\Omega,\calf,\P, W)$ for the optimal switching
problem one can always obtain a setting
for a randomized optimal control problem by the following direct construction.
Take an arbitrary probability space
$(\Omega',\calf',\P')$ where
  a  Poisson random
 measure $\mu$ with intensity $\lambda$ is defined.
 Thus in particular,
 for every $\omega'\in\Omega'$, $\mu(\omega', dt\,da)$
 is a measure on $(0,\infty)\times A$.
 Let us define $\hat \Omega=\Omega\times \Omega'$, let us denote by $\hat \calf$ the completion of
the product $\sigma$-algebra
$\calf\otimes \calf'$ with respect to $\P\otimes \P'$ and by
$\hat \P$ the extension of $\P\otimes \P'$ to $\hat \calf$.
One can introduce canonical extensions
$\hat W$ and $\hat \mu$
of $W$ and $\mu$ to $\hat\Omega$  by setting
$$\hat W_t(\omega,\omega')=W_t(\omega),
\qquad \hat\mu(\omega,\omega', dt\,da)=\mu(\omega', dt\,da),
$$ for every $t\ge 0$, $\omega\in\Omega$,
$\omega'\in\Omega'$.  Then it can be easily   checked
that,  under $\hat \P$, $\hat W$ is a standard Wiener process
and $\hat\mu$
 is a random Poisson measure on $(0,\infty)\times A$
 with the same intensity $\lambda$,
 independent of $\hat W$.  So we see
 that $(\hat \Omega, \hat \calf,\hat \P, \hat  W, \hat  \mu)$
is a setting for a randomized control problem as formulated before, that we call
\emph{product extension} of the setting $(\Omega,\Fc,\P,W)$ for the initial optimal
switching problem. This construction will be used again
for the proofs of several results below and will be further studied.
\newline
We note that by a classical result,
 see for instance \cite{Zabczyk96} Theorem 2.3.1, we may take $\Omega'=[0,1]$,
 $\calf'$ the corresponding Borel sets and $\P'$ the Lebesgue measure. This shows
that the extended setting is rather ``economical'' in the loose sense
that it does not introduce much randomness with respect to the original setting.
\newline
We also note that the initial formulation of a randomized
setting 
was more general, since  it was not required that
$\hat \Omega$ should be a product space $\Omega\times \Omega'$ and, even if it were the case, it was not required that the process $ \hat W $
should depend only on $\omega\in\Omega$ while the random measure $\hat\mu$ should depend only on  $\omega'\in\Omega'$.
}
\ep
\end{Remark}

\subsection{Equivalence of the optimal switching  and the randomized control pro\-blems}
We can now state one of the main results of the paper.

\begin{Theorem}
\label{MainThm}
Assume that {\bf (A1)} and {\bf (A2)} are satisfied.
Then the values of the optimal switching problem
and of the randomized control problem are equal:
\beq \label{equivrandom}
{\text{\Large$\upsilon$}}_0 &=& {\text{\Large$\upsilon$}}_0^\Rc,
\enq
where ${\text{\Large$\upsilon$}}_0$ and ${\text{\Large$\upsilon$}}_0^\Rc$ are defined by  \eqref{primalvalue} and \eqref{dualvalue} respectively.
This common value only depends on the objects $A,b,\sigma,f,g,c,x_0,\xi_0$
appearing in assumption {\bf (A1)}.
\end{Theorem}

The last sentence follows immediately from
Remark
\ref{indepofthesetting}, from the equality ${\text{\Large$\upsilon$}}_0$ $=$  ${\text{\Large$\upsilon$}}_0^\Rc$
and from the obvious fact that ${\text{\Large$\upsilon$}}_0$ cannot depend on
the measure $\lambda$ introduced in assumption {\bf (A2)}.
The following section is entirely devoted to the proof of the equality.

\section{Proof of Theorem \ref{MainThm}}
\label{proofthm}

\subsection{Preliminaries}
\label{prelimproofthm}

In this section, {\bf (A1)} and {\bf (A2)} are always assumed to hold.
We will prove separately the two inequalities ${\text{\Large$\upsilon$}}_0^\Rc$ $\le$ ${\text{\Large$\upsilon$}}_0$ and ${\text{\Large$\upsilon$}}_0$ $\le$ ${\text{\Large$\upsilon$}}_0^\Rc$. In both cases, we need
similar constructions, which consist in starting with a given
setting $(\Omega,\calf,\P,  W)$ for the  optimal switching problem
 formulated in section \ref{Primal},
  building a product space  by adding another suitable
  probability space as an independent factor and thus arriving
at a suitable setting  $(\hat \Omega, \hat \calf,\hat \P,  \hat  W, \hat  \mu)$
 for a  randomized control problem as formulated before. In this paragraph
we present this construction and its main properties needed later.

\bigskip

Let us start with a setting
 $(\Omega,\calf,\P, W)$, where $(\Omega,\calf,\P)$ is a complete probability
 space and $W$ a $d$-dimensional standard Wiener process and let $(\Omega',\calf',\P')$ be another arbitrary probability space.
We finally set $\hat \Omega=\Omega\times \Omega'$ and denote by $\hat \calf$ the completion of
the product $\sigma$-algebra
$\calf\otimes \calf'$ with respect to $\P\otimes \P'$ and by
$\hat \P$ the extension of $\P\otimes \P'$ to $\hat \calf$.
One can introduce a canonical extension of $W$ to $\hat\Omega$ setting
$\hat W_t(\omega,\omega')=W_t(\omega)$ for every $t\ge 0$, $\omega\in\Omega$,
$\omega'\in\Omega'$.  Then $\hat W$ is a standard Wiener process under $\hat \P$,
as it can be easily checked.
More generally, any random element defined in $\Omega$ or $\Omega'$ has
an extension defined by similar formulae, whose  law under $\hat \P$ is the
same as the law under the original probability.

\bigskip

One can formulate an optimal switching problem in the new setting
$(\hat \Omega, \hat \calf,\hat \P, \hat  W)$ in the same way as before:
we let $ \F^{\hat W}=( \calf^{\hat W}_t)_{t\ge 0}$ denote the right-continuous and
$\hat \P$-complete filtration generated by $\hat W$,
and we define the set of admissible strategies $\hat\cala$
as the elements of the form
 $\hat\alpha=(\hat\tau_n,\hat\xi_n)_{n\ge 1}$
satisfying properties analogous to $(i)-(v)$ in section \ref{Primal}, but
with the filtration $ \F^{\hat W}$ instead of $ \F^{ W}$. For
 any
$\hat\alpha\in \hat\cala$ one
finds the corresponding trajectory $\hat X^{\hat \alpha}$ solving
 the controlled equation
\begin{equation}\label{stateeqext}
    d\hat X_t^{\hat \alpha} \ = \ b_t( \hat X^{\hat \alpha}, \hat\alpha(t))\,dt +
\sigma_t( \hat X^{\hat \alpha}, \hat \alpha(t))\,d\hat W_t,
\qquad
 \hat X_0^{\hat \alpha}=x_0,
\end{equation}
where $\hat \alpha(\cdot)$ is  the piecewise constant process
associated to $\hat \alpha$,
and computes
the corresponding reward:
\begin{equation}\label{rewardext}
\hat J(\hat \alpha)  :=   \hat\E\Big[\int_0^Tf_t(\hat X^{\hat \alpha},\hat\alpha(t))\,dt
+g(\hat X^{\hat \alpha},\hat \alpha(T))\Big]-
\hat  \E\Big[\sum_{n\ge 1}1_{\hat \tau_n<T}\,
c_{\hat \tau_n}(\hat X^{\hat \alpha},\hat \xi_{n-1},\hat \xi_n)\Big].
\end{equation}
Finally, the value is defined as
\begin{equation}\label{primalvalueext}
{\text{\Large$\hat \upsilon$}}_0 := \sup_{\hat\alpha\in\hat\Ac} \hat J(\hat\alpha).
\end{equation}

\bigskip

One may wish to compare this value with the value of the switching problem
formulated in the original setting  $(\Omega,\calf,\P, W)$.
To this end,
let us recall that  $\F^W=(\calf^W_t)_{t\ge 0}$ denotes the right-continuous and
$\P$-complete filtration in $\Omega$ generated by $W$. Every $\sigma$-algebra $\calf_t^W$
gives rise to a $\sigma$-algebra in $\hat \Omega$ defined as
$$
\calf_t^W \times \Omega':=\{A\times \Omega'\;:\; A\in\calf_t^W\}.
$$
This way one obtains a new filtration in $\hat \Omega$ (which is right-continuous
but not $\hat\P$-complete in general).  Recalling that
$  \F^{\hat W}=(  \calf^{\hat W}_t)_{t\ge 0}$ denotes the right-continuous and
$\hat \P$-complete filtration generated by $\hat W$, and letting
$\caln$ denote the family of $\hat \P$-null sets in $\hat \calf$, one arrives at the equality
\begin{equation}\label{descrizfiltrW}
    \calf^{\hat W}_t= (\calf_t^W \times \Omega')\vee \caln, \qquad t\ge0,
\end{equation}
which can be verified by lengthy but standard arguments.

If $\tau$ is an $\F^W$-stopping time then its canonical extension
defined by $\hat\tau(\omega,\omega')=\tau(\omega)$ is
a $\F^{\hat W}$-stopping time; indeed, for every $t\ge 0$,
$\{\hat \tau\le t\}=\{\tau\le t\}\times \Omega'$ belongs
to $\calf_t^W \times \Omega'$ and so to $ \calf^{\hat W}_t$.
Now suppose that
  $A\in\calf^W_\tau$; then
for every $t\ge 0$,
$$
(A\times \Omega')\cap \{\hat \tau\le t\}=
(A\cap
\{\tau\le t\})\times \Omega'\in \calf_t^W \times \Omega'\subset
 \calf^{\hat W}_t.
$$
This shows that if  $A\in\calf^W_\tau$ then
 $A\times \Omega'\in  \calf^{\hat W}_{\hat \tau}$.
 This property implies that for any $\calf^W_\tau$-measurable random
 variable $\xi$, its canonical extension $\hat\xi(\omega,\omega')=\xi(\omega)$
 is $ \calf^{\hat W}_{\hat \tau}$-measurable.

It follows that if we start from an admissible control strategy
$\alpha\in\Ac$  of the form
$
\alpha=(\tau_n,\xi_n)_{n\ge 1}$ and denote
$\hat\tau_n,\hat\xi_n$ the canonical extensions of
$\tau_n,\xi_n$ respectively, then
$
\hat\alpha:=(\hat\tau_n,\hat\xi_n)_{n\ge 1}$
is an admissible strategy for the optimal switching problem formulated
in the setting $(\hat\Omega,\hat \calf,\hat P,\hat W)$, hence an element
of $\hat \cala$.
Moreover, it is easy to realize that in this case the process
$\hat X^{\hat \alpha}$ solution to  \eqref{stateeqext} is the same as
  the canonical extension of the process $X^\alpha$
defined  as the solution
to the controlled equation \eqref{stateeq}
($\hat X^{\hat \alpha}_t(\omega,\omega')= X_t^\alpha(\omega)$) and,
moreover,
 the reward
\eqref{rewardext} is the same as the original one: $\hat J(\hat \alpha)= J( \alpha)$.
We deduce that the two values satisfy the inequality
$$
{\text{\Large$\upsilon$}}_0 = \sup_{\alpha\in\Ac} J(\alpha)
\le  \sup_{\hat\alpha\in\hat \Ac} \hat J(\hat\alpha)=
{\text{\Large$\hat \upsilon$}}_0.
$$

\bigskip

Following \cite{BCFP16bAAP}, we next introduce a variant of the optimal switching problem
formulated in the new setting
$(\hat \Omega, \hat \calf,\hat \P, \hat  W)$. We define a new filtration,
denoted $ \F^{\hat W,\infty}=(\calf^{\hat W,\infty}_t)_{t\ge 0}$,
as follows: we first introduce
$$
\Omega\times \calf' :=\{\Omega \times B\;:\; B\in\calf'\},
$$
which is a  $\sigma$-algebra in $\hat \Omega$, and then set
$$
 \calf^{\hat W,\infty}_t:=
 \calf^{\hat W}_t \vee (\Omega\times \calf')= (\calf_t^W \times \Omega')\vee \caln
\vee (\Omega\times \calf'), \qquad t\ge0.
$$
Next we define a new set of admissible strategies, denoted $\hat\cala^\infty$,
consisting of
  the elements of the form
 $\hat\alpha=(\hat\tau_n,\hat\xi_n)_{n\ge 1}$
satisfying properties analogous to $(i)-(v)$ in section \ref{Primal}, but
with the filtration $ \F^{\hat W,\infty}$ instead of $ \F^{ W}$. For any such
$\hat\alpha$ one
finds the corresponding trajectory $\hat X^{\hat \alpha}$ solving
 the controlled equation
 \eqref{stateeqext}  and computes
the corresponding reward $\hat J(\hat\alpha)$ by
 \eqref{rewardext} as before.
 The corresponding value is defined as
\begin{equation}\label{primalvalueextinfty}
{\text{\Large$\hat \upsilon$}}_0^\infty := \sup_{\hat\alpha\in\hat\Ac^\infty} \hat J(\hat\alpha).
\end{equation}
Since $ \F^{\hat W}$ is a smaller filtration than
$ \F^{\hat W,\infty}$,
we have $\hat\Ac \subset \hat\Ac^\infty $ and we conclude that
$
{\text{\Large$\upsilon$}}_0 \le
{\text{\Large$\hat \upsilon$}}_0 \le
{\text{\Large$\hat \upsilon$}}_0^\infty.
$
Actually, it turns out that the three values in fact coincide:

\begin{Lemma}\label{L:AWmu}
With the previous notations we have
$
{\text{\Large$\upsilon$}}_0 ={\text{\Large$\hat \upsilon$}}_0 =
{\text{\Large$\hat \upsilon$}}_0^\infty.
$
 \end{Lemma}

The intuitive explanation is that in the optimal switching problem
for ${\text{\Large$\hat \upsilon$}}_0^\infty$ the controller has
access to the information coming from the Wiener filtration
as well as the one represented by the $\sigma$-algebra $\Omega\times\calf'$;
however, under $\hat\P$ the latter is independent of $W$
and so it has no use in getting a better performance.
We do not write down the proof of this Lemma, since it is entirely
analogous to Lemma 4.1 of
\cite{BCFP16bAAP} (there the notation  $  \F^{ W,\mu'_\infty}$
and $ \Ac^{W,\mu'} $
was used instead of our notation $ \F^{\hat W,\infty}$ and $\hat\Ac^\infty $).
The conclusion of this lemma will be used
in the proof of the inequality
${\text{\Large$\upsilon$}}_0^\Rc$ $\le$ ${\text{\Large$\upsilon$}}_0$
below.

\subsection{Proof of the inequality ${\text{\Large$\upsilon$}}_0^\Rc$ $\le$ ${\text{\Large$\upsilon$}}_0$}\label{firstineq}

We follow closely \cite{BCFP16bAAP},
making use in particular of the basic  Proposition 4.2
in that paper.

Let $(\Omega,\calf,\P, W)$
be a setting for the  optimal switching problem  formulated in section \ref{Primal}.
We construct a setting for a randomized control problem in the form of an
appropriate product
extension as described in Remark   \ref{productextension}.

Let $\lambda$ be a Borel measure on $A$ satisfying {\bf (A2)}. As a first step, we
construct a suitable surjective measurable map $\pi:\R\to A$ and
 a  measure $\lambda'$ on the Borel subsets of the real line
satisfying the condition  $\lambda$ $=$ $\lambda' \circ \pi^{-1}$
(the image measure of $\lambda'$ under $\pi$) and such that
$\lambda'(\{r\})=0$ for every $r\in\R$.

We do not report the details of the construction of $\pi$ and $\lambda'$,
for which we refer the reader to paragraph 4.1 of \cite{BCFP16bAAP}. We just
mention that it is a very simple consequence of the well known fact that
the space of modes $A$, being a Borel space, is  known to be either finite or countable (with the discrete topology) or isomorphic, as a measurable space, to the real line: see e.g. \cite{BertsekasShreve78}, Corollary 7.16.1.

Next, we choose  $(\Omega',\calf',\P')$ to be the canonical probability space of a
 non-explosive Poisson point process on $(0,\infty)\times \R$ with intensity $\lambda'$.
 Thus, $\Omega'$ is the set of sequences $\omega'=(t_n,r_n)_{n\geq1}\subset(0,\infty)\times \R$ with $t_n<t_{n+1}\nearrow\infty$,
 $(\sigma_n,\rho_n)_{n\geq1}$ is the canonical marked point process
 (i.e.  $\sigma_n(\omega')=t_n$, $\rho_n(\omega')=r_n$),  and
 $\mu'$ $=$ $\sum_{n\ge 1}\delta_{(\sigma_n,\rho_n)}$ is the corresponding random measure.
 Let $\Fc'$ denote the  smallest $\sigma$-algebra such that all the maps
 $\sigma_n,\rho_n$ are measurable, and   $\P'$  the unique probability on
$\Fc'$ such that  $\mu'$ is a Poisson random measure with intensity $\lambda'$ (since $\lambda'$ is a finite measure,  this probability actually exists).
 We will also use the completion of the space  $(\Omega',\calf',\P')$, still denoted by the same symbol by abuse of notation. Setting
$$
\eta_n=\pi(\rho_n), \qquad
\mu=\sum_{n\ge 1}\delta_{(\sigma_n,\eta_n)},
$$
it is easy to verify that $\mu$ is a Poisson random measure on $(0,\infty)\times A$ with intensity $\lambda$, defined in $(\Omega',\calf',\P')$.

Then we perform the construction
described in section \ref{prelimproofthm}: we
define $\hat \Omega=\Omega\times \Omega'$, we denote by $\hat \calf$ the completion of
$\calf\otimes \calf'$ with respect to
$\P\otimes \P'$ and by
$\hat \P$ the extension of $\P\otimes \P'$ to $\hat \calf$.
As explained before, $W$ has a canonical extension to
a $\hat\P$-standard Wiener process $\hat W$ in
$\hat \Omega$. The Poisson random measure $\mu$ also has
a canonical extension to a random measure $\hat\mu$ on $(0,\infty)\times A$
defined on $\hat\Omega$ setting $\hat\mu= \sum_{n\ge 1}\delta_{(\hat\sigma_n,\hat\eta_n)}$,
where 
$\hat\sigma_n(\omega,\omega'):= \sigma_n(\omega')$  and $ \hat\eta_n(\omega,\omega'):=\eta_n(\omega')$.
It is immediate to verify that $\hat\mu$ is also a Poisson
random measure with intensity $\lambda$, independent of $\hat W$.
We may summarize this construction saying that
$(\hat \Omega, \hat \calf,\hat \P,  \hat  W, \hat  \mu)$
is a setting for a  randomized control problem.

We can then formulate the corresponding randomized optimization problem
 as in section \ref{randomizedformulation}:
we define the   $\hat\P$-completed filtration
$\F^{\hat W,\hat \mu}=(\calf^{\hat W,\hat \mu}_t)_{t\ge 0}$
generated by $\hat W$ and $\hat \mu$
as in formula
\eqref{filtraznaturaleWmu},
we   introduce the classes $  \hat\calv, \hat\Vc_{\inf\,>\,0}$ and, for any
 admissible control $ \hat\nu \in \hat\calv$,
 the corresponding martingale $\kappa^{\hat\nu}$, the probability
 $\hat\P^{\hat \nu}(d\omega\, d\omega' )=\kappa_T^{\hat\nu}(\omega,\omega')\,\hat\P(d \omega\,d\omega')$,
 the processes $\hat I$ and $\hat X$,
 given by formula
\eqref{I} and solution to  \eqref{dynXrandom}
respectively,
  the reward $ J^\Rc(\hat\nu)$ given by
 \eqref{defJrandomized}-\eqref{JRuno}-\eqref{JRdue}
and the value
${\text{\Large$\upsilon$}}_0^\Rc $
defined in \eqref{dualvalue}.
We recall that this value does not depend on the specific setting
chosen above for the randomized optimal control problem,
as noticed in Remark \ref{indepofthesetting}.

We mention that we have the following alternative description of the
filtration  $\F^{\hat W,\hat \mu}\!=\!(\calf^{\hat W,\hat \mu}_t)_{t\ge 0}$.
We first
introduce in $(\Omega',\Fc')$ the $\P'$-complete right-continuous
  filtration $\F^\mu=(\calf^\mu_t)_{t\ge 0}$,
generated by $\mu$ and defined  by
$$
\calf^\mu_t = \sigma\big(\mu((0,s]\times C)\,:\, s\in [0,t],\, C\in\calb(A)\big)
\vee \caln',
$$
where $\caln'$ denotes the family of $\P'$-null sets of $\Fc'$.
Next we
introduce the $\sigma$-algebra in $\hat \Omega$ defined as
$$
\Omega \times \calf^\mu_t:=\{ \Omega \times B\;:\; B\in\calf^\mu_t\}.
$$
Then we have the equality
\begin{equation}\label{filtraznaturalealternativa}
        \calf^{\hat W,\hat\mu}_t= (\calf_t^W \times \Omega')\vee (\Omega \times \calf^\mu_t)
    \vee \caln, \qquad t\ge0,
\end{equation}
which is analogous to formula
\eqref{descrizfiltrW} and can be proved by similar arguments.

\bigskip

At this point we make use of the following technical result,
which is a special case of Proposition 4.2 in \cite{BCFP16bAAP}:

\begin{Proposition}
\label{P:Ineq_I}
For every $\hat\nu\in\hat\Vc_{\inf\,>\,0}$ there exists
$\hat\alpha^{\hat\nu}\in\Ac^{\hat W,\infty}$ such that
\begin{equation}\label{LawB_I=LawB_alpha}
\mathscr L_{\hat\P^{\hat \nu}}(\hat W,\hat I) \ =
\ \mathscr L_{\hat\P}(\hat W,\hat\alpha^{\hat\nu}),
\end{equation}
i.e., the law of $(\hat W,\hat I)$ under $\hat\P^{\hat \nu}$ is the same as the
law of $(\hat W,\hat\alpha^{\hat\nu})$ under $\hat\P$.
\end{Proposition}

The proof of the inequality
${\text{\Large$\upsilon$}}_0^\Rc$ $\le$  ${\text{\Large$\upsilon$}}_0$ is now
finished as follows. Take $\hat\nu\in \hat\Vc_{\inf\,>\,0}$ and construct
$\hat\alpha^{\hat\nu}\in\Ac^{\hat W,\infty}$ as in Proposition
\ref{P:Ineq_I}.
Since $\hat X$ is obtained solving equation \eqref{dynXrandom} and
  $\hat X^{\hat\alpha^{\hat\nu}}$ is obtained solving equation \eqref{stateeqext}
  (with $\hat\alpha^{\hat\nu}$ instead of $\hat\alpha$)
  it is a well-known fact that under the conditions in
 Assumption {\bf (A1)}  the equality  \eqref{LawB_I=LawB_alpha}
 implies that
\begin{equation}
\mathscr L_{\hat \P^{\hat \nu}}(\hat X,\hat I)
=  \mathscr L_{\hat\P}(\hat X^{\hat\alpha^\nu},\hat\alpha^{\hat\nu}).
\end{equation}
This immediately entails that
$J^\Rc(\hat \nu)=\hat J(\hat\alpha^{\hat \nu})$. It follows that
$J^\Rc(\hat \nu)\le {\text{\Large$\hat \upsilon$}}_0^\infty$, where
the latter was defined in \eqref{primalvalueextinfty}.
 From the arbitrariness of $\hat \nu$ we deduce that $
 \sup_{\hat\nu\in\hat\Vc_{\inf\,>\,0}}J^\Rc(\hat\nu)
 \leq {\text{\Large$\hat \upsilon$}}_0^\infty$.
 From  \eqref{eqinfzero}  it follows that
  ${\text{\Large$\upsilon$}}_0^\Rc \leq{\text{\Large$\hat \upsilon$}}_0^\infty$.
  Since by
Lemma \ref{L:AWmu} we have $
{\text{\Large$\upsilon$}}_0 =
{\text{\Large$\hat \upsilon$}}_0^\infty
$
we arrive at the desired conclusion
  ${\text{\Large$\upsilon$}}_0^\Rc$ $\leq $ ${\text{\Large$\upsilon$}}_0$.
\ep

\subsection{Proof of the inequality ${\text{\Large$\upsilon$}}_0$ $\le$ ${\text{\Large$\upsilon$}}_0^\Rc$}
\label{secVleqVR}

In this proof we borrow some constructions from
 \cite{FP15} and \cite{BCFP16bAAP}, but
 the proofs need substantial extensions.
 This is due to the cost $J_2(\alpha)$ 
 in \eqref{Jdue} 
 related to switching: in the various approximations and convergence
 arguments in Lemmas  \ref{approxuno}, 
  \ref{approxdue},  \ref{approxtre}, 
  this term requires a special treatment, in particular
  because we do not require boundedness
  of the costs $c_t(x,a,a')$ 
  nor a uniform bound on the number $N_T$ of 
  switchings. 
Additional difficulties are also due to the fact that the terminal reward $g(X^\alpha, \alpha(T))$  depends on the final mode $\alpha(T)$.
 

Suppose we are given a setting $(\Omega,\calf,\P, \F, W)$ for the optimal switching problem
as described in section \ref{Primal}, and  consider the controlled equation \eqref{stateeq} and the
reward \eqref{gaineq}.

\begin{Lemma}\label{approxuno}
For any $\delta >0$ there exists an admissible switching strategy
$\alpha=(\tau_n,\xi_n)_{n\ge 1}\in\cala$
such that
$$
J(\alpha)\ge {\text{\Large$\upsilon$}}_0-\delta
$$
and moreover
\begin{enumerate}
\item[(i)] there exists an integer $N\ge 1$ such that $\tau_n =+\infty$ as soon as $n >N$,
\item[(ii)] the set $\{\xi_n(\omega)\; :\; \omega\in\Omega, n=1,\ldots, N\}$ is finite.
\end{enumerate}
 \end{Lemma}

For the proof of this Lemma, we need the following stability result,
that will be used several times below.
Following \cite{80Krylov},
for any pair  $\alpha^1,\alpha^2:\Omega\times [0,T]\to A$ of measurable
  processes we define a distance $\tilde\rho(\alpha^1,\alpha^2)$ setting
\beqs
\tilde \rho(\alpha^1,\alpha^2) &=& \E \Big[\int_0^T\rho(\alpha^1_t,\alpha^2_t)\,dt \Big].
\enqs
where $\rho$ is an arbitrary metric compatible with the topology of $A$
and satisfying $\rho<1$.
Using in particular
the continuity condition {\bf (A1)}-(iii) one can show the following.

\begin{Lemma}
\label{contrhotilde}
Suppose we have a probability space $(\hat\Omega,
\hat \calf,\Q)$ with filtrations
$\G^k=(\calg_t^k)_{t\ge 0}$ ($k\ge 0$)
and a process $B$ which is a Wiener process with respect to each $\G^k$.
Consider
the equations
$$
d  Y_t ^k=  b_t( Y^k, \gamma^{k}(t))\,dt
+ \sigma_t(Y^k, \gamma^{k}(t))\,dB_t,
\qquad
\hat Y_0^k =x_0,
$$
where each $\gamma^k$ is an admissible switching strategy
with respect to $\G^k$ (i.e., satisfying
 properties  $(i)-(v)$ in section \ref{Primal}, but
with the filtration $\G^k$ instead of $ \F^{ W}$).
Suppose that
 \begin{equation}\label{perconvstab}
 \tilde\rho (\gamma^k, \gamma^0)\to 0,
 \quad{\rm and}\quad
 \gamma^{k}(T)\to \gamma^{0}(T) \;\Q-a.s.
 \end{equation}
as $k\to\infty$.
Then for every $p\in [1,\infty)$,
\begin{equation}\label{primaconvstab}
    \E^\Q\sup_{t\in [0,T]}|Y^k_t-Y^0_t|^p\to 0,
\quad
\E^\Q \Big[ \int_0^Tf_t(Y^k,\gamma^{k}(t))\,dt \Big]
\to
\E^\Q \Big[ \int_0^Tf_t(Y^0,\gamma^{0}(t))\,dt \Big].
\end{equation}
\begin{equation}\label{secondaconvstab}
\E^\Q \Big[  g(Y^k,\gamma^{k}(T))\Big]
\to
\E^\Q \Big[  g(Y^0,\gamma^{0}(T))\Big],
\end{equation}
so that in particular $J_1(\gamma^k)\to J_1(\gamma^0)$.
\end{Lemma}
\noindent {\bf Proof.}
The convergence  result \eqref{primaconvstab}
 was first proved in \cite{80Krylov} in the standard diffusion case.
The simple extension to the non-Markovian case is presented
 in \cite{FP15}, Lemma 4.1 and Remark 4.1. This holds under
 the condition $\tilde\rho (\gamma^k, \gamma^0)\to 0$ alone.
 Using the second assumption in
  \eqref{perconvstab}, the continuity assumption
 {\bf (A1)}-(iii) and the growth conditions \eqref{PolGrowth_f_g},
 the convergence
  \eqref{secondaconvstab} follows easily.
\ep

\bigskip

\noindent {\bf Proof of Lemma \ref{approxuno}.}
By the definition of ${\text{\Large$\upsilon$}}_0$, for any $\delta >0$
there exists an admissible switching strategy
$\alpha=(\tau_n,\xi_n)_{n\ge 1}\in\cala$
such that
$
J(\alpha)\ge {\text{\Large$\upsilon$}}_0-\delta/3.
$
Next we modify $\alpha$ in two steps, in order to satisfy the additional
requirements in the statement of the Lemma.

In a first step we consider the strategy obtained by taking only the
first $N$ switchings in $\alpha$, that we denote
 $\alpha^N=(\tau_n,\xi_n)_{n= 1}^N$. Formally, we use
 this notation to indicate the strategy where we have modified
the pairs $(\tau_n,\xi_n)$ for $ n>N$ setting them equal to $(  \infty,\bar \xi)$
where $\bar\xi\in A$ is fixed arbitrarily. We claim that $J(\alpha^N)\ge J(\alpha)-2\delta/3$
for $N$ sufficiently large.

To verify the claim we first note that, for the piecewise constant processes $\alpha^N(\cdot)$,
$\alpha(\cdot)$
associated to $\alpha^N$ and $\alpha$ we have
$\alpha^N(t)$ $=$
$\alpha(t)$ for $t\in [0,T\wedge \tau_N]$ and so
$$
\tilde \rho(\alpha^N(\cdot),\alpha(\cdot)) =
\E \Big[\int_0^T\rho(\alpha^N(t),\alpha(t))\,dt \Big]
=
\E \Big[\int_{\tau_N\wedge T}^T\rho(\alpha^N(t),\alpha(t))\,dt \Big]
\le
\E \,[T- ({\tau_N\wedge T})]\to 0,
$$
since $\tau_N\to\infty$. Since $\{\tau_N>T\}$ $\subset$
$\{\alpha^N(T)= \alpha(T)\}$ we also have
$\P(\alpha^N(T)= \alpha(T))\geq \P(\tau_N>T)\to 1$,
so that $\alpha^N(T)\to  \alpha(T)$ in $\P$-probability
and, passing to a subsequence if necessary, we may assume
$\alpha^N(T)\to  \alpha(T)$ $\P$-a.s.
Applying Lemma \ref{contrhotilde} to the controlled equations
satisfied by $  X^{ \alpha^N}$ and
$X^{ \alpha}$ 
and setting $B=  W$, $Y^k=  X^{ \alpha^k}$,
$\gamma^k(\cdot)=  \alpha^k(\cdot)$
and $Y^0=  X^{  \alpha}$, $\gamma^0(\cdot)=  \alpha(\cdot)$, 
we conclude that
$  J_1( \alpha^N)\to   J_1(  \alpha)$.

Since $\alpha^N(t)$ $=$
$\alpha(t)$ for $t\in [0,T\wedge \tau_N]$ we also have
$X^{\alpha^N}_t $ $=$
$X^\alpha_t$ for $t\in [0,T\wedge \tau_N]$ and
therefore for $n=1,\ldots,N$ we have
$$
1_{\tau_n<T}\,c_{\tau_n}(X^{\alpha^N},\xi_{n-1},\xi_n)=
1_{\tau_n<T}\,c_{\tau_n}(X^\alpha,\xi_{n-1},\xi_n).
$$
If $N$ is chosen so large that   $| J_1( \alpha^N) -   J_1(  \alpha)|<\delta/3$
then, taking into account the fact that costs are nonnegative, we obtain
\beqs
J(\alpha)=J_1(\alpha)-J_2(\alpha)
&\le &
J_1(\alpha) -
  \E\Big[\sum_{n= 1}^N1_{\tau_n<T}\,c_{\tau_n}(X^\alpha,\xi_{n-1},\xi_n)\Big]
  \\
&= &
J_1(\alpha) -
  \E\Big[\sum_{n= 1}^N1_{\tau_n<T}\,c_{\tau_n}(X^{\alpha^N},\xi_{n-1},\xi_n)\Big]
\\
&\le &
  J_1(\alpha^N)+\delta/3 -
  \E\Big[\sum_{n= 1}^N1_{\tau_n<T}\,c_{\tau_n}(X^{\alpha^N},\xi_{n-1},\xi_n)\Big]
\\
&=&
  J (\alpha^N)+\delta/3,
\enqs
and since we have $
J(\alpha)\ge {\text{\Large$\upsilon$}}_0-\delta/3
$ we obtain
$
J(\alpha^N)\ge {\text{\Large$\upsilon$}}_0-2\delta/3
$ as claimed.

\bigskip

As a second  step we fix $N$ and we further modify $\alpha^N$
in the following way.  Since $A$ is a Borel space, it is separable.
Let us fix a dense sequence $(a_i)_{i\ge 1}$ and define, for each
integer $k\ge 1$, a map $\Pi_k :A\to A$ that
assigns to each $b\in A$ its nearest point in $\{a_1,\ldots,a_k\}$,
more precisely
$$
\Pi_k(b)=a_{i(b)}, \quad
{\rm where}\quad i(b):=\min \{
j\in \{1,\ldots,k\} \; :\; \rho(b,a_j)\le \rho(b,a_i) {\rm \; for \; all\;}
i\in\{ 1,\ldots, k\}\}.
$$
It is easy to see that $\Pi_k :A\to A$ is Borel measurable and
$\rho (\Pi_k(a),a)\downarrow 0$ as $k\to\infty$.

Starting from the strategy $\alpha^N=(\tau_n,\xi_n)_{n= 1}^N$ constructed
above we define $\alpha^{N,k}=(\tau_n,\Pi_k(\xi_n))_{n= 1}^N$.
We note that each strategy $\alpha^{N,k}$ satisfies the conditions
stated in the Lemma. To finish the proof it is therefore enough
to prove that $J(\alpha^{N,k})\to J(\alpha^{N})$ as $k\to\infty$: indeed, taking any
$k$ sufficiently large we have $J(\alpha^{N,k})\ge J(\alpha^{N})-\delta/3$
so that any such strategy satisfies
$J(\alpha^{N,k})\ge  {\text{\Large$\upsilon$}}_0-\delta$.

In order to prove that $J(\alpha^{N,k})\to J(\alpha^{N})$
we start noting that $\rho (\alpha^{N,k}(t),\alpha^{N}(t))\to 0$ $\P$-a.s.
for every $t\in[0,T]$. In particular $\alpha^{N,k}(T)\to\alpha^{N}(T)$ $\P$-a.s.
and we also have  $\tilde\rho (\alpha^{N,k}(\cdot),\alpha^{N}(\cdot))\to 0$.
Another application of Lemma \ref{contrhotilde} shows that
$J_1(\alpha^{N,k})\to J_1(\alpha^{N})$ and we also have
$$ \forall \; p\in [1,\infty), \quad \quad   \E \sup_{t\in [0,T]}|X^{\alpha^{N,k}}_t-X^{\alpha^{N}}_t|^p\to 0.$$ 
Passing to a subsequence if necessary, we may assume that
$\sup_{t\in [0,T]}|X^{\alpha^{N,k}}_t-X^{\alpha^{N}}_t|\to 0$ $\P$-a.s. and for
every $n=1,\ldots,N$ we have, by the continuity assumptions in
{\bf (A1)}-(iii),
$$
1_{\tau_n<T}\,c_{\tau_n}(X^{\alpha^{N,k}},\Pi_k(\xi_{n-1}),\Pi_k(\xi_n))
\to
1_{\tau_n<T}\,c_{\tau_n}(X^{\alpha^N},\xi_{n-1},\xi_n),
\qquad \P-a.s.
$$
From the growth condition \eqref{PolGrowth_f_g} we obtain the inequality
$$
0\le 1_{\tau_n<T}\,c_{\tau_n}(X^{\alpha^{N,k}},\Pi_k(\xi_{n-1}),\Pi_k(\xi_n))
\le 1_{\tau_n<T}\, L\, (1+ \sup_{t\in [0,T]}|X^{\alpha^{N,k}}|)^r
$$
and by \eqref{growthsolprimal} we conclude that the right-hand side is bounded
in $L^p(\P)$ for every $p\in[1,\infty)$. It follows that
$$
\E\Big[1_{\tau_n<T}\,c_{\tau_n}(X^{\alpha^{N,k}},\Pi_k(\xi_{n-1}),\Pi_k(\xi_n))\Big]
\to
\E\Big[1_{\tau_n<T}\,c_{\tau_n}(X^{\alpha^N},\xi_{n-1},\xi_n)\Big],
$$
and we conclude that $J_2(\alpha^{N,k})\to J_2(\alpha^{N})$ since the number
of switchings is bounded by $N$. This way we have proved that
$J(\alpha^{N,k})\to J(\alpha^{N})$, which ends the proof of the Lemma.
\ep

\bigskip

In order to proceed further we need to construct a product probability space
as explained in section \ref{prelimproofthm},
making use of a properly
chosen auxiliary
probability space denoted
 $(\Omega',\calf',\P')$. This can be taken as an arbitrary  probability space
 where appropriate random objects  are defined.
For integers $m,n,k\ge1$,  we assume that real random variables  $ U_n^m $, $S^m_n$
and random measures
 $\pi^k$ are defined on  $(\Omega',\calf',\P')$ and satisfy the following conditions:
\begin{enumerate}
\item every $ U_n^m $ is uniformly  distributed on $(0,1)$;

\item  every   $ S_n^m$ admits
 a density (denoted $f_n^m(t)$)  with respect to the Lebesgue measure,
 and we have $ 0<S_1^m<S_2^m<S^m_3<\ldots $ for every $m$, and 
 $ S_n^m\to 0$ as $m\to\infty$
 for every $n$;

\item  every $\pi^k$ is a Poisson random measure on $(0,\infty)\times A$, admitting compensator $k^{-1}\lambda(da)\,dt$ with respect to its natural filtration;
 \item the random elements  $U^m_n,S^h_j$, $\pi^k$    are all independent.
 \end{enumerate}
The inequalities required in point 2. above can be satisfied for instance by
choosing the support of each density $f_n^m$ inside the interval
$((1-2^{-n})/m, (1-2^{-n-1})/m)$.
The role of these random elements will become clear in the constructions that
follow. Notice that for the construction of the space
$(\Omega',\calf',\P')$ only the knowledge of the measure $\lambda$ is required.
Moreover  by a classical result,
 see \cite{Zabczyk96} Theorem 2.3.1, we may take $\Omega'=[0,1]$,
 $\calf'$ the corresponding Borel sets and $\P'$ the Lebesgue measure.

Next we perform the construction described in section \ref{prelimproofthm}.
Let us define $\hat \Omega=\Omega\times \Omega'$, let us denote by $\hat \calf$ the completion of
the product $\sigma$-algebra
$\calf\otimes \calf'$ with respect to $\P\otimes \P'$ and by $\Q$
 the extension of $\P\otimes \P'$ to $\hat \calf$ (the notation $\hat \P$ will be
used for a different probability introduced below).
 As before we denote
$\hat W_t$,  $\hat U^m_n$, $\hat S^h_j$, $\hat \pi^k$
 the canonical extensions of $W$, $U^m_n,S^h_j$, $\pi^k$
to $\hat\Omega$.

 Since $\hat W$ is a standard Wiener process under $\hat \P$ we can consider the
optimal switching problem in the setting $(\hat \Omega, \hat \calf, \hat \P,\hat W)$
as in section   \ref{prelimproofthm}:
we define the set of admissible strategies $\hat\cala$
as the elements of the form
 $\hat\alpha=(\hat\tau_n,\hat\xi_n)_{n\ge 1}$
satisfying properties analogous to $(i)-(v)$ in section \ref{Primal}, but
with the filtration $ \F^{\hat W}$ instead of $ \F^{ W}$.  For
 any
$\hat\alpha\in \hat\cala$ one
finds the corresponding trajectory $\hat X^{\hat \alpha}$ solving
 the controlled equation
 \eqref{stateeqext}
and computes
the corresponding reward $\hat J(\hat \alpha)  $ given in
 \eqref{rewardext}, namely
\begin{equation}\label{rewardextdue}
\hat J(\hat \alpha)  = \hat J_1(\hat \alpha) -\hat J_2(\hat \alpha)
=
  \E^\Q\Big[\int_0^Tf_t(\hat X^{\hat \alpha},\hat\alpha(t))\,dt
+g(\hat X^{\hat \alpha},\hat \alpha(T))\Big]-
  \E^\Q\Big[\sum_{n\ge 1}1_{\hat \tau_n<T}\,
c_{\hat \tau_n}(\hat X^{\hat \alpha},\hat \xi_{n-1},\hat \xi_n)\Big],
\end{equation}
where $\E^\Q$ denotes the expectation under $\Q$.
It was explained in Remark \ref{strategiecomemisure} that
any switching strategy can be viewed as a marked point process in $A$.
In the following it will be convenient to identify any
$\hat\alpha\in\hat\cala$
of the form
 $\hat\alpha=(\hat\tau_n,\hat\xi_n)_{n\ge 1}$  with the corresponding
random measure on $(0,\infty)\times A$ defined as
$$
\hat\alpha=\sum_{n\ge 1}\delta _{(\hat\tau_n,\hat\xi_n)}\,1_{\hat\tau_n<\infty}
$$
where $\delta$ denotes the Dirac measure. We will use the same symbol
to denote the strategy and the corresponding measure.
We will also need the corresponding natural filtration
$\F^{\hat \alpha}=(\calf^{\hat \alpha}_t)_{t\ge 0}$ in $(\hat\Omega,\hat\calf)$
defined  by the formula:
\beq\label{filtraznaturalealpha}
\calf^{ \hat \alpha}_t&=&\sigma (
\hat \alpha((0,s]\times C)\,:\, s\in [0,t],\, C\in\calb(A)),
\enq
and also the filtration
$\F^{\hat W}\vee \F^{\hat\alpha}:=(\calf_t^{\hat W}\vee\calf_t^{\hat\alpha})_{t\geq 0}$.
We denote $\calp(\F^{\hat \alpha})$, $\calp(\F^{\hat W}\vee \F^{\hat\alpha})$,
the corresponding predictable
$\sigma$-algebras.

A basic role in the arguments below will be played by the concept of
compensator (or dual predictable projection) of this random measure,
as presented for instance in \cite{ja}.

\begin{Lemma}\label{approxdue}
For any $\delta >0$ there exists an admissible switching strategy
$\hat{{\beta}}\in\hat\cala$
such that
$$
\hat J(\hat \beta)\ge {\text{\Large$\upsilon$}}_0-2\delta
$$
and moreover the $\Q$-compensator of the corresponding random measure
on $(0,T]\times A$
with respect to $\F^{\hat W}\vee \F^{\hat\beta}$ is absolutely continuous
with respect to the measure $\lambda (da)\,dt$ and it has the form
$$
\hat\nu_t^{\hat\beta}(\omega,\omega',a)\,\lambda (da)\,dt
$$
where $\hat\nu^{\hat\beta}: \hat\Omega\times [0,T]\times A\to [0,\infty)$
 is a  $\calp(\F^{\hat W}\vee \F^{\hat\beta})\otimes \calb(A)$-measurable
 function.
 \end{Lemma}

{\bf Proof.}
 Given $\delta >0$, let us consider the strategy $\alpha$ constructed in
 Lemma \ref{approxuno} and let us denote $\hat\alpha= (\hat\tau_n,\hat\xi_n)_{n\ge 1}$
 its canonical extension. We have seen in section \ref{prelimproofthm} that
 $\hat\alpha\in\hat\cala$ and  $J(\alpha)=\hat J(\hat\alpha)$.
 By construction of $\hat \alpha$ it holds that
$\hat J(\hat \alpha)\ge {\text{\Large$\upsilon$}}_0-\delta$,
$\hat\tau_{n} =\infty$ as soon as $n >N$, and
 the set $\{\hat \xi_n(\omega)\; :\; \omega\in\Omega, n=1,\ldots, N\}$ is finite. The
 corresponding random measure and piecewise constant process are
$$\hat\alpha=\sum_{n= 1}^N\delta _{(\hat\tau_n,\hat\xi_n)}\,1_{\hat\tau_n<\infty},
\qquad
\hat\alpha(t) \; = \; \xi_01_{ [0,\hat \tau_{1})}(t)
+\sum_{n=1}^N
\hat \xi_{n}1_{ [\hat \tau_n,\hat \tau_{n+1})}(t),
$$
where $\xi_0\in A$ is the given starting mode.

The idea of the proof is to perturb this random measure slightly
in such a way that the corresponding reward will not be changed too much
and at the same time its compensator will have the desired properties.

Let $\rho$ be a metric inducing the topology of $A$ and satisfying $\rho < 1$.
For every   $m\ge 1$, let $\bfB(b,1/m)$ denote the open ball  of radius $1/m$, with respect to the metric $\rho$, centered at $b\in A$.
Since $\lambda(da)$ has full support, we have $\lambda(\bfB(b,1/m))>0$ and we can  define a transition kernel $q^m(b,da)$ in $A$ setting
\beqs
q^m(b,da)&=& \frac{1}{\lambda(\bfB(b,1/m))}\, 1_{\bfB(b,1/m)}(a) \lambda(da).
\enqs
We recall that we require $A$ to be a Borel space,
and we denote by $\calb(A)$ its Borel $\sigma$-algebra.
There exists a Borel measurable function $q^m:A\times [0,1]\to A$
such that for every $b\in A$ the measure $B\mapsto q^m(b,B)$ ($B\in \calb(A)$) is the image of the Lebesgue measure  on $[0,1]$ under the mapping
$u\mapsto q^m(b,u)$.
Thus, if  $U$ is a random variable defined on some probability space and having uniform law on $[0,1]$ then, for fixed $b\in A$, the $A$-valued random variable
$q^m(b,U)$ has law $q^m(b,da)$. The use of the same symbol $q^m$ should not generate confusion.
The existence of the function $q^m$ (even for a general transition kernel on $A$)
 is well known when $A$ is a
 separable complete metric space,
in particular, when $A$ is the unit interval $[0,1]$,
(see e.g. \cite{Zabczyk96}, Theorem 3.1.1) and the general case reduces to this one,
since it is known that any Borel  space is either finite or countable (with the discrete topology)
or isomorphic, as a measurable space, to the interval  $[0,1]$: see e.g.
\cite{BertsekasShreve78}, Corollary 7.16.1.

\bigskip

For fixed   $m\ge 1$,
define    $\hat R^m_0=0$ and
 $$
\hat R_n^m=\hat\tau_n+\hat S^m_n,
\quad \hat \beta_n^m=q^m(\hat \xi_n,\hat U^m_n),
\qquad \quad  n\ge 1.
$$
Since   $\hat \tau_n<\hat\tau_{n+1}$ and since $\hat S^m_n>0$  
we see that $\hat\alpha^m:=(\hat R^m_n,\hat\beta^m_n)_{n\ge1}$ is an
admissible strategy (the property that $\Q(\hat R^m_n=T$ for some $n)=0$
comes from the fact that $\hat S_n^m$ have absolutely continuous
laws and are independent of $\hat\tau_n$).
Let
$$\hat\alpha^m=\sum_{n=1}^N\delta_{(\hat R^m_n,\hat \beta^m_n)},
\qquad
\hat\alpha^m(t) \; = \; \xi_01_{ [0,\hat R^m_{1})}(t)
+\sum_{n=1}^N
\hat \beta^m_{n}1_{ [\hat R^m_n,\hat R^m_{n+1})}(t),
$$
denote the corresponding random measure and the associated
piecewise constant process.

It is possible to compute explicitly the
$\Q$-compensator of these random measures
with respect to $\F^{\hat W}\vee \F^{\hat\alpha^m}$, which
is given by the formula
$$
  \sum_{n=1}^N  1_{(\hat\tau_n\vee \hat R^m_{n-1},\hat R^m_n]}(t)\,
 q^m(\hat\xi_n,da) \frac{f_n^m(t-\hat \tau_n)}{1-F_n^m(t-\hat \tau_n)}\,dt,
$$
where we denote  by  $F_n^m(s)=\int_{-\infty}^sf_n^m(t)dt$ the cumulative distribution function of $S_n^m$, with  the convention that $\frac{f_n^m(s)}{1-F^m_n(s)}=0$ if $F^m_n(s)=1$.
The proof of this result is given in
 Lemma A.11 in \cite{FP15}.
We can write this formula in the form
  $$
\left[\sum_{n= 1}^N  1_{(\hat\tau_n\vee \hat R^m_{n-1},\hat R^m_n]}(t)\,
\frac{1}{\lambda(\bfB(\hat \xi_n,1/m))}\, 1_{\bfB(\hat \xi_n,1/m)}(a)
\frac{f_n^m(t-\hat \tau_n)}{1-F_n^m(t-\hat \tau_n)}
\right]
\, \lambda(da)\,dt
 $$
where the function in square brackets is
a nonnegative $\calp(\F^{\hat W}\vee \F^{\hat\alpha^m})\otimes \calb(A)$-measurable  function.

To finish the proof it is enough to show that $\hat J(\hat\alpha^m)\to \hat J(\hat\alpha)$
as $m\to\infty$ (or at least for a subsequence $m_k$).
Indeed, since $\hat J(\hat \alpha)\ge {\text{\Large$\upsilon$}}_0-\delta$,
 for large $m$ we will have
$\hat J(\hat\alpha^m) \ge{\text{\Large$\upsilon$}}_0-2\delta$
and we can take $\hat\beta= \hat\alpha^m$ for such $m$ in the statement
of the Lemma, since its compensator satisfies the required conditions.

To prove the required convergence
$\hat J(\hat\alpha^m)\to \hat J(\hat\alpha)$
we first note that 
$$
0<\hat R_n^m -\hat \tau_n= \hat S^m_n\to 0,
\qquad \Q-a.s.
$$
We deduce that $\Q$-a.s., $\hat\alpha^m(t)\to \hat\alpha(t) $,
except perhaps at points $\hat\tau_n$, and so $dt$-a.s.
In particular, since there are no switchings at the terminal time $T$,
we have $\Q(\hat\tau_n=T {\rm\;for\;some\;} n)=0$ and we conclude
that
$\hat\alpha^m(T)\to \hat\alpha(T) $  $\Q$-a.s.
We also note that by the choice of the kernel $q^m(b,da)$
we have $\rho(\hat \xi_n,\hat \beta^m_n)<1/m\to 0$ and
therefore
 for the distance already considered above we have
\begin{equation}\label{compensatorofmodifiedppmbis}
\tilde\rho (\hat{{\alpha}}, \hat\alpha^m)=
\E^\Q \Big[\int_0^T\rho(\hat{{\alpha}}(t),\hat{{\alpha}}^m(t))\,dt \Big]
\to 0, \qquad m\to\infty.
\end{equation}
Applying Lemma \ref{contrhotilde} to the controlled equations
satisfied by $\hat X^{\hat\alpha^m}$ and
$\hat X^{\hat\alpha}$ 
 and setting 
$B=\hat W$, $Y^k=\hat X^{\hat \alpha^k}$,
$\gamma^k(\cdot)=\hat \alpha^k(\cdot)$
and $Y^0=\hat X^{\hat \alpha}$, $\gamma^0(\cdot)=\hat \alpha(\cdot)$
we conclude that
$\hat J_1(\hat \alpha^m)\to \hat J_1(\hat \alpha)$.

It remains to study the convergence of
$\hat J_2(\hat \alpha^m)$. Since it is a finite sum, it is enough to check that
for every $n=1,\ldots,N$
\begin{equation}\label{c_nconvergono}
    \E^\Q\Big[ 1_{\hat R^m_n<T}\,
c_{\hat R^m_n}(\hat X^{\hat \alpha^m},\hat \beta^m_{n-1},\hat \beta^m_n)\Big]\to
   \E^\Q\Big[ 1_{\hat \tau_n<T}\,
c_{\hat \tau_n}(\hat X^{\hat \alpha},\hat \xi_{n-1},\hat \xi_n)\Big],
\end{equation}
as $m\to\infty$. By the growth condition \eqref{PolGrowth_f_g} in {\bf (A1)} we have
$$
|c_{\hat R^m_n}(\hat X^{\hat \alpha^m},\hat \beta^m_{n-1},\hat \beta^m_n)|
\le
L \,(1+ \sup_{t\in [0,T]}|X_t^{\hat \alpha^m}|)^r
$$
and the right-hand side is bounded in all $L^p(\Q)$ spaces, by the estimate
\eqref{growthsolprimal}. So it is enough to check that we have
convergence
 $\Q$-almost surely  
for the terms in right brackets
in \eqref{c_nconvergono}.
 Once again, since $\tilde\rho (\hat{{\alpha}}, \hat\alpha^m) \to 0$, the application of  Lemma \ref{contrhotilde} gives that, for any $p \in [1, \infty)$
$$
\E^\Q[\sup_{t\in [0,T]}|X_t^{\hat \alpha^m}- X_t^{\hat \alpha}|^p]\to 0,$$
and so, at least
for a subsequence, we have  $\|X^{\hat \alpha^m}- X^{\hat \alpha}\|_\infty\to 0$
$\Q$-a.s.
We have already checked above that, $\Q$-a.s.,
$\hat R^m_n\to \hat\tau_n$,  $\hat \beta^m_n\to \hat \xi_n$
and so we have
$1_{\hat R^m_n<T}\to 1_{\hat \tau_n<T}$
and finally
$$
1_{\hat R^m_n<T}\,
c_{\hat R^m_n}(\hat X^{\hat \alpha^m},\hat \beta^m_{n-1},\hat \beta^m_n)\to
1_{\hat \tau_n<T}\,
c_{\hat \tau_n}(\hat X^{\hat \alpha},\hat \xi_{n-1},\hat \xi_n),
$$
by the continuity properties of the coefficient $c$ stated in Assumption {\bf (A1)}.
The required convergence \eqref{c_nconvergono} is proved and the proof of
the Lemma is finished.
\ep

\begin{Lemma}\label{approxtre}
For any $\delta >0$ there exists an admissible switching strategy
$\hat\alpha\in\hat\cala$
such that
$$
\hat J(\hat \alpha)\ge {\text{\Large$\upsilon$}}_0-3\delta
$$
and moreover the $\Q$-compensator of the corresponding random measure
on $(0,T]\times A$
with respect to $\F^{\hat W}\vee \F^{\hat\alpha}$ is absolutely continuous with respect to the measure $\lambda (da)\,dt$ and it has the form
$$
\hat\nu_t(\omega,\omega',a)\,\lambda (da)\,dt
$$
where $\nu: \hat\Omega\times [0,T]\times A\to [0,\infty)$
 is  a $\calp(\F^{\hat W}\vee \F^{\hat\alpha})\otimes \calb(A)$-measurable
 function satisfying
 $
 \inf\nu >0.
 $
 Moreover, denoting by $N_T $ the number of jump times
of $\hat\alpha$
in $[0,T]$, we have $N_T\in L^p(\Q) $
for
every $p\in [1,\infty)$.
 \end{Lemma}

{\bf Proof.} Let
$\hat{{\beta}}\in\hat\cala$ be the switching strategy constructed in
 Lemma \ref{approxdue}, that we write in the form of a random measure
$\hat \beta=\sum_{n= 1}^N1_{\hat R_n<\infty}\delta_{(\hat R_n,\hat \beta_n)}$
having at most $N$ summands.
 The idea of the proof is to modify  the associated random
measure   by adding an independent Poisson process
with ``small'' intensity. This will not affect  the
reward too much and will
produce a random measure whose compensator
remains absolutely continuous with respect to the measure
$\lambda(da)\,dt$ with a bounded density which,  in addition, is bounded
away from zero.

Recall that on the space
$(\Omega',\calf',\P')$ we assumed that for every integer $k\ge 1$ there
existed    a Poisson random measure $\pi^k$  on $(0,\infty)\times A$, admitting compensator $k^{-1}\lambda(da)\,dt$ with respect to its natural filtration.
We denoted  $\hat\pi^k$  its canonical extension to $(\hat\Omega,\hat\calf)$,
that we
write
in the form of a random measure on $(0,\infty)\times A$:
\beqs
\hat \pi^k &=& \sum_{n\ge 1}\delta_{(\hat\sigma_n^k , {\hat \eta_n^k })},
\enqs
for a marked point process $(\hat\sigma_n^k , {\hat\eta_n^k })_{n\ge1}$.
Let us define other random measures setting
$$
\hat\mu^{k} =\hat\beta+\hat\pi^k.
$$
Note that the jump times $(\hat R_n)_{n\ge1}$  are independent
of the jump times $(\hat\sigma_n^k )_{n\ge1}$, and the latter have
absolutely continuous laws. It follows that, except possibly
on a set of $\Q$ probability zero, their graphs  are
disjoint, i.e. $\hat\beta$ and $\hat\pi^k $ have no common jumps,
and $\hat\mu^k$ do not charge the terminal time $T$. Therefore,
  the random measures $\hat \mu^{k} $ can be identified
with admissible switching strategies (they belong to $\hat\cala$)
and, together with their 
associated piecewise constant processes
(denoted $\hat  \mu^{k}(\cdot)$)
    admit a representation of the form
$$
\hat\mu^{k} =\sum_{n\ge 1}\delta_{(\hat \tau_n^{k} , {\hat\xi^{k} _n})},
\qquad
\hat \mu^{k}(t) = \xi_0 1_{ [0 ,\hat\tau^{k}_{1})}(t)+
 \sum_{n\ge 1 } \xi^{k} _{n}1_{ [\hat\tau_n^{k} ,\hat\tau^{k} _{n+1})}(t), \;\;\; t \in[ 0,T],
$$
where $\xi _{0}$ is the starting mode,
$(\hat\tau_n^{k} ,\hat \xi^{k} _n)_{n\ge 1}$ is a marked point process,
each $\hat\tau_n^{k} $ coincides with one of the times $\hat R_n$
or one of the times
$\hat\sigma_n^k $, and each $\hat\xi_n^{k} $
coincides with one of the random variables $\hat\eta_n^k $ or one of the random variables
$\hat \beta_{n}$.

We recall that $\hat\beta$ had at most $N$  switchings, and we define
$N^k_T := \sum_{n\ge 1}1_{\hat\sigma_n^k\le T}$ which
has Poisson law with parameter $\lambda(A)T/k$. It follows that the number
of jump times $\hat\tau_n^{k}$ in $[0,T]$ of each
$\hat\mu^{k}$ cannot exceed $N+N_T^k$
and therefore it belongs to $ L^p(\Q) $
for
every $p\in [1,\infty)$.

Let us verify that the
$\Q$-compensator of each $\hat\mu^k$
with respect to $\F^{\hat W}\vee \F^{\hat\mu^k}$ satisfies
the properties in the statement of the Lemma.
We first note that,
since $\hat\beta$ and $\hat\pi^k $ are independent, it is easy to prove
 that
 $\hat\mu^k=\hat \beta+\hat \pi^{k} $
  has compensator  $(\hat\nu_t^{\hat\beta}(\omega,\omega',a)+k^{-1})\,\lambda(da)\,dt$
with respect  to the filtration
$\F^{\hat W}\vee \F^{\hat\beta}\vee  \F^{\hat \pi^k}$=$
(\calf_t^{\hat W}\vee \calf^{\hat\beta}_t\vee \calf^{\hat\pi^k }_t)_{t\ge0}$.
Let us denote  $\F^{\hat \pi^k}$ $=$ $(\calf^{\hat\pi^k}_t)$ the natural filtration
of $\hat\pi^k$
defined as in \eqref{filtraznaturalealpha}.
We wish to compute the $\Q$-compensator of $\hat\mu^k$ with respect to
the filtration
$\F^{\hat W}\vee \F^{\hat\mu^k}$
$=$ $( \calf_t^{\hat W}\vee\calf_t^{\hat\mu^k})_{t\geq 0}$, which is smaller
than $\F^{\hat W}\vee \F^{\hat\beta}\vee  \F^{\hat \pi^k}$.
To this end, consider the measure space
$([0,\infty)\times \hat \Omega \times A, \calb([0,\infty))\otimes \hat \calf \otimes \calb(A),dt\otimes\Q(d\hat \omega)\otimes\lambda(da))$.
Although this is not a probability space, one can define in a standard way the conditional expectation
 of any positive measurable
function, given an arbitrary sub-$\sigma$-algebra.
Let us denote by  $\hat\nu_t^k(\omega,\omega',a)$ the conditional expectation
of the  random field $ \hat\nu_t^{\hat\beta}(\omega,\omega',a)+k^{-1}$
 with respect to the $\sigma$-algebra $\calp(\F^{\hat W}\vee \F^{\hat\mu^k})
 \otimes \calb(A)$. It is then easy to verify that
the compensator of $\hat\mu ^k$ with respect to
$\F^{\hat W}\vee \F^{\hat\mu^k}$ coincides with $\hat\nu^k$.
Moreover, since  $ \hat\nu_t^{\hat\beta}$ is nonnegative,
we can take a version of $\hat\nu^k$ satisfying
$$
\inf_{\hat\Omega\times [0,T]\times A}\hat\nu^k \ge  k^{-1} >0.
  $$

To finish the proof of Lemma
\ref{approxtre} it is enough to show that $\hat J(\hat\mu^k)\to \hat J(\hat\beta)$
as $k\to\infty$ (or at least for a subsequence).
Indeed, since $\hat J(\hat\beta) \ge{\text{\Large$\upsilon$}}_0-2\delta$,
for large $k$ we will have
$\hat J(\hat\mu^k) \ge{\text{\Large$\upsilon$}}_0-3\delta$
and we can take $\hat\alpha= \hat\mu^k$ for such $k$ in the statement
of the Lemma, since its compensator satisfies the required conditions.

 We first claim that, for large $k$, $\hat\mu^{k}(\cdot)$
is close to $\hat\beta(\cdot)$ with respect to the metric $\tilde \rho$,
namely that
\begin{equation}\label{approxvalue}
\tilde\rho(\hat\mu^{k}(\cdot),\hat\beta(\cdot)):=
\E^\Q \Big[\int_0^T\rho(\hat\mu^{k}(t),\hat\beta(t))\,dt \Big]
\to 0, \qquad k\to\infty.
\end{equation}
Recall that the jump times of $\hat\pi^k $ are denoted $\hat\sigma_n^k $.
Since  $\hat\sigma_1^k $ has exponential law with parameter
$ \lambda(A)/k$ the event $B_k =\{\hat\sigma_1^k   >T\}$ has probability
$e^{-\lambda(A)T/k}$, so that $\Q(B_k)\to 1$ as $k\to \infty$.
We note that, on the set $B_k$, we have $\hat\mu^k(t) =
\hat\beta(t)$ for all $t\in [0,T]$. Since   we assume $\rho <1$,
we have  $\tilde\rho(\mu^{k}(\cdot),\hat\beta(\cdot))\le T (1-\Q(B_k))$ and
 the claim
 \eqref{approxvalue} follows immediately.

Similarly, since  $\hat\mu^k(T)=\hat\beta(T)$ on $B_k$, we have
$\hat\mu^k(T)\to\hat\beta(T)$ in $\Q$-probability, and passing to a subsequence
(denoted by the same symbol)
if necessary we can assume $\hat\mu^k(T)\to\hat\beta(T)$ $\Q$-a.s.

 Applying Lemma \ref{contrhotilde} to the controlled equations
satisfied by $\hat X^{\hat\mu^k}$ and
$\hat X^{\hat\beta}$ 
and setting  
$B=\hat W$, $Y^k=\hat X^{\hat \alpha^k}$,
$\gamma^k(\cdot)=\hat \alpha^k(\cdot)$
and $Y^0=\hat X^{\hat \alpha}$, $\gamma^0(\cdot)=\hat \alpha(\cdot)$
we conclude that
$\hat J_1(\hat \mu^k)\to \hat J_1(\hat \beta)$.
It remains to study the convergence of
$$
\hat J_2(\hat \mu^k)=\E^\Q\Big[\sum_{n\ge 1}1_{\hat \tau_n^k<T}\,
c_{\hat \tau^k_n}(\hat X^{\hat \mu^k},\hat \xi^k_{n-1},\hat \xi^k_n)\Big].
$$
We recall that $\hat\beta$ had at most $N$ switchings, and we defined
$N^k_T  = \sum_{n\ge 1}1_{\hat\sigma_n^k\le T}$ which
has Poisson law with parameter $\lambda(A)T/k$.
  By the growth conditions in Assumption {\bf (A1)} we have
$$
  \sum_{n\ge 1}1_{\hat \tau_n^k<T}\,
c_{\hat \tau^k_n}(\hat X^{\hat \mu^k},\hat \xi^k_{n-1},\hat \xi^k_n)
\le (N+N^k_T)\, L\, (1+\sup_{t\in [0,T]}|\hat X^{\hat\mu^k}  _t|)^r
  $$
and recalling \eqref{growthsolprimal}   we see that for
every $p\in [1,\infty)$ the right-hand side
is bounded in $L^p(\Q)$ by a constant independent of $k$.
Setting again $B_k =\{\hat\sigma_1^k  >T \}$ and recalling that
  $\Q(B_k)\to 1$, by the H\"{o}lder inequality we conclude that
$$  \E^\Q\Big[1_{B_k^c}\;\sum_{n\ge 1}1_{\hat \tau_n^k<T}\,
c_{\hat \tau^k_n}(\hat X^{\hat \mu^k},\hat \xi^k_{n-1},\hat \xi^k_n)\Big]
\to 0, \qquad k\to\infty.
$$
Next we note that on the event $B_k$ the measures $\hat \mu^k$ and $\hat\beta$
coincide on $(0,T]\times A$ and therefore on $B_k\times [0,T]$
we also have $\hat\mu^k(\cdot) =
\hat\beta(\cdot)$ and  $\hat X^{\hat \mu^k}=\hat X^{\hat \beta}$ $\Q$-a.s.
It follows that
\beqs
\hat J_2(\hat \mu^k)&=&\E^\Q\Big[1_{B_k}\;\sum_{n\ge 1}1_{\hat \tau_n^k<T}\,
c_{\hat \tau^k_n}(\hat X^{\hat \mu^k},\hat \xi^k_{n-1},\hat \xi^k_n)\Big]
+
\E^\Q\Big[1_{B_k^c}\;\sum_{n\ge 1}1_{\hat \tau_n^k<T}\,
c_{\hat \tau^k_n}(\hat X^{\hat \mu^k},\hat \xi^k_{n-1},\hat \xi^k_n)\Big]
\\
&=&\E^\Q\Big[1_{B_k}\;\sum_{n= 1}^N1_{\hat R_n<T}\,
c_{\hat R_n}(\hat X^{\hat\beta},\hat \beta_{n-1},\hat \beta_n)\Big]
+
\E^\Q\Big[1_{B_k^c}\;\sum_{n\ge 1}1_{\hat \tau_n^k<T}\,
c_{\hat \tau^k_n}(\hat X^{\hat \mu^k},\hat \xi^k_{n-1},\hat \xi^k_n)\Big]
\\
&\le &\hat J_2(\hat \beta)
+
\E^\Q\Big[1_{B_k^c}\;\sum_{n\ge 1}1_{\hat \tau_n^k<T}\,
c_{\hat \tau^k_n}(\hat X^{\hat \mu^k},\hat \xi^k_{n-1},\hat \xi^k_n)\Big].
\enqs
Since we clearly have $\hat J_2(\hat \beta)\le \hat J_2(\hat \mu^k)$
it follows that $ \hat J_2(\hat \mu^k)\to \hat J_2(\hat \beta)$.  Now we
have verified that $\hat J(\hat\mu^k)\to \hat J(\hat\beta)$
and the proof of Lemma
\ref{approxtre}
is finished.
\ep

\bigskip

We are now able to end the proof of the inequality
${\text{\Large$\upsilon$}}_0$ $\le$ ${\text{\Large$\upsilon$}}_0^\Rc$.

Let $\delta>0$ be given and denote $\hat\alpha=\sum_{n\ge 1}\delta _{(\hat\sigma_n,\hat\eta_n)}\,1_{\hat\sigma_n<\infty}$
the random measure corresponding  to the strategy  $\hat\alpha$ given by Lemma
\ref{approxtre}.

Let
$\caln$ denote the family of $\Q$-null sets of $(\hat\Omega,\hat\calf)$.
Then the filtration
$(\calf^{\hat W}_t\vee\calf_t^{\hat\alpha}\vee \caln)_{t\geq 0}$
coincides with  the filtration previously denoted by
$\F^{\hat W,\hat \alpha}=(\calf^{\hat W,\hat \alpha}_t)_{t\ge 0}$
(compare with formula  \eqref{filtraznaturaleWmu} or \eqref{filtraznaturalealternativa}).
It is easy to see that
$ \hat\nu_t ( \omega,\omega',a)\, \lambda(da)\,dt$ is the $\Q$-compensator
of $\hat\alpha$  with respect to $\F^{\hat W,\hat \alpha}$ as well.

Using the Girsanov theorem for point processes (see e.g. \cite{ja}) we next
construct an equivalent probability under which $\hat\alpha$   becomes
a Poisson random measure with intensity $\lambda$.
Since the function $\hat\nu$ occurring in  Lemma \ref{approxtre} is a strictly positive
$\calp(\F^{\hat W,\hat\alpha})\otimes\calb(A)$-measurable
random field with bounded inverse, the Dol\'eans exponential process
\begin{equation}\label{defgirsanovdensity}
 M_t \; := \; \exp\Big(\int_0^t\int_A(1-\hat\nu_s(a)^{-1})\,
 \hat\nu_s(a) \lambda(da)\,ds\Big) \prod_{\hat \sigma_n\le t} \hat\nu_{\hat \sigma_n}
 (\hat\eta_n)^{-1},\qquad t\in [0,T],
\end{equation}
is a strictly positive martingale (with respect to $\F^{\hat W,\hat\alpha}$ and $\Q$),
and we can
define an equivalent probability $\hat\P$ on the space $(\hat\Omega,\hat\calf)$ setting
 $\hat\P(d\omega\,d\omega')$ $=$ $M_T(\omega,\omega')\Q(d\omega\,d\omega')$.
 The expectation under
 $\hat\P$ will be denoted  $\hat\E$. By the Girsanov theorem,
 the restriction of
 $\hat\alpha$      to $(0,T]\times A$
  has $(\hat\P,\F^{\hat W,\hat\alpha})$-compensator $\lambda(da)\,dt$,
so that  in particular it is a Poisson random measure.
It can also be proved by standard arguments (see e.g. \cite{FP15},
page 2155,
for   detailed verifications in a similar framework)  that
$\hat W$ remains a $(\hat\P,\F^{\hat W,\hat\alpha})$-Wiener process and
that $\hat W$ and $\hat \alpha$ are independent under $\hat\P$.
We have thus constructed a setting
$(\hat \Omega, \hat \calf,\hat \P,  \hat W, \hat \alpha)$ for a randomized
control problem as in section \ref{randomizedformulation}.


Although the random field
$\nu$ is not bounded in general, so in particular it does not
belong
to the class $\hat \calv$ of admissible controls for the
randomized control problem,
  we can still
introduce the
 Dol\'eans exponential process $\kappa^{\hat \nu}$
 corresponding to $\hat\nu$ by the formula
\eqref{doleans}, namely:
\begin{equation}\label{doleansbis}
\kappa_t^{\hat\nu} =
 \exp\left(\int_0^t\int_A (1 - \hat \nu_s(a))\lambda(da)\,ds
\right)\prod_{\hat \sigma_n \le t} \hat \nu_{\hat \sigma_n}(\hat \eta_n),\qquad t\in [0,T].
\end{equation}
 Comparing \eqref{defgirsanovdensity}  and
 \eqref{doleansbis}  shows that
  $\kappa^{\hat \nu}_T\,M_T\equiv 1$.
It follows that
$\hat\E [\kappa^{\hat \nu}_T]=\E^\Q [M_T\kappa^{\hat \nu}_T]=1$, so that
$\kappa^{\hat \nu}$ is indeed
a  $\hat\P$-martingale on $[0,T]$ and   we can
define the corresponding probability
$\hat\P^{\hat\nu}(d\hat \omega):=\kappa^{\hat \nu}_T(\hat \omega)\hat\P(d\hat \omega)$.
Since
  $\kappa^{\hat \nu}_T\,M_T\equiv 1$,
 the  Girsanov transformation $\hat\P\mapsto \hat\P^{\hat\nu}$
is the inverse of the transformation $\Q\mapsto \hat\P$ made above,
and changes back the probability $\hat\P$ into
$\Q$ considered above, so that we have
$\hat\P^{\hat\nu}=\Q$.

Let       $\hat X$ be the
solution to
the equation
\begin{equation}\label{controlledhatkappa}
d\hat X_t=  b_t( \hat X, \hat I_t)\,dt
+ \sigma_t(\hat X, \hat I_t)\,d \hat W_t,
\qquad
\hat X_0 = x_0,
\end{equation}
where $\hat I$ is the piecewise constant $A$-valued process associated to $\hat\alpha$
and starting at the initial mode $\xi_0$ (the same as in formula \eqref{I},
and elsewhere indicated $\hat\alpha(\cdot)$):
\begin{equation}\label{Idue}
\hat I_t \ =
\xi_0\,1_{[0,\hat \sigma_{ 1})}(t)+
\sum_{n\ge 1}\hat \eta_n\,1_{[\hat \sigma_n,\hat \sigma_{n+1})}(t), \qquad t\ge 0.
\end{equation}
The corresponding reward of the switching problem is then
\beq\nonumber
\hat J(\hat \alpha) &=&  \E^\Q \left[
 \int_0^Tf_t(\hat X,\hat I_t)\,dt+g(\hat  X,\hat I_T)
-
\sum_{n\ge 1}1_{\hat \sigma_n<T}\,c_{\hat \sigma_n}(\hat X,\hat \eta_{n-1},\hat \eta_n)
 \right]
 \\\label{JRtre}
 &=&
 \hat\E^{\hat\nu}
\left[\int_0^Tf_t(\hat X,\hat I_t)\,dt+g(\hat  X,\hat I_T)
-
\sum_{n\ge 1}1_{\hat \sigma_n<T}\,c_{\hat \sigma_n}(\hat X,\hat \eta_{n-1},\hat \eta_n)
\right],
\enq
where we have used  $\hat\P^{\hat\nu}=\Q$ in the last equality.

For any integer $k\ge 1$, let define $\hat\nu^k_t(a)=\hat\nu_t(a)\wedge k$. Therefore
$\hat\nu^k\in\hat\calv$,  we can define the corresponding process
$\kappa^{\hat \nu^k}$
 by formula \eqref{doleans},
 the probability
 $\hat\P^{\hat\nu^k}(d\hat \omega)=\kappa^{\hat \nu^k}_T(\hat \omega)\,\hat\P(d\hat \omega)$, and
compute the reward $J^\Rc(\hat\nu^k)$ of the corresponding  randomized
problem.  Since equation \eqref{controlledhatkappa} coincides with
the randomized equation \eqref{dynXrandom}, this is given by
\beq \label{JRk}
J^\Rc(\hat\nu^k) &=&  \hat\E^{\hat\nu^k}
\left[\int_0^Tf_t(\hat X,\hat I_t)\,dt+g(\hat  X,\hat I_T)
-
\sum_{n\ge 1}1_{\hat \sigma_n<T}\,c_{\hat \sigma_n}(\hat X,\hat \eta_{n-1},\hat \eta_n)
\right],
\enq
where $\hat\E^{\hat\nu^k} $ denotes the expectation under $\hat\P^{\nu^k}$.

We claim that $J^\Rc(\hat\nu^k) \to \hat J(\hat\alpha)$ as $k\to\infty$.
Assuming this for a moment,
since
$\hat J(\hat \alpha)\ge {\text{\Large$\upsilon$}}_0-3\delta$,
we will have
 $J^\Rc(\hat\nu^k)  \ge {\text{\Large$\upsilon$}}_0-4\delta$ for large $k$, and
since $J^\Rc(\hat\nu^k)$ is the reward of a randomized control problem,
by Remark \ref{indepofthesetting} it can not exceed the value
${\text{\Large$\upsilon$}}_0^\Rc$ defined in \eqref{dualvalue}, whatever
the setting where the randomized problem is formulated.
It follows that
$
   {\text{\Large$\upsilon$}}_0^\Rc \ge  {\text{\Large$\upsilon$}}_0-4\delta
$
and by the arbitrariness of $\delta $  we obtain the required inequality
 ${\text{\Large$\upsilon$}}_0^\Rc$ $\ge$ ${\text{\Large$\upsilon$}}_0$.

It remains to prove the
claim that $J^\Rc(\hat\nu^k) \to \hat J(\hat\alpha)$.
Setting
$$
\Phi= \int_0^Tf_t(\hat X,\hat I_t)\,dt+g(\hat  X,\hat I_T)
-
\sum_{n\ge 1}1_{\hat \sigma_n<T}\,c_{\hat \sigma_n}(\hat X,\hat \eta_{n-1},\hat \eta_n)
$$
and
comparing
\eqref{JRtre} with \eqref{JRk}, proving the claim amounts to showing
that $\hat\E^{\hat\nu^k}[\Phi]\to \hat\E^{\hat\nu}[\Phi]$ or equivalently
$\hat\E[\kappa^{\hat \nu^k}_T\Phi]\to \hat\E [\kappa^{\hat \nu}_T\Phi]$.
Using the growth condition \eqref{PolGrowth_f_g} in Assumption {\bf (A1)}
we see that
$$
| \Phi|
\le c\, (1+N_T)\,  (1+\sup_{t\in [0,T]}|\hat X_t|)^r
  $$
  for a suitable constant $c$,
  where
$N_T $ denotes the number of jump times $\hat\sigma_n$
of $\hat\alpha$
in $[0,T]$. From Lemma \ref{approxtre}
we know that  $N_T\in L^p(\Q) $
for
every $p\in [1,\infty)$ and by \eqref{EstimateX_nu}
we conclude that
$\Phi\in L^p(\Q) $
for
every $p\in [1,\infty)$ as well. The required convergence
$\hat\E[\kappa^{\hat \nu^k}_T\Phi]\to \hat\E [\kappa^{\hat \nu}_T\Phi]$
can now be verified by standard arguments, exactly the same
as in   \cite{FP15},
pages 2156-2157.

\ep

\section{The randomized BSDE}
\label{Sec:separandom}

In this section the assumptions {\bf (A1)} and {\bf (A2)}
are assumed to hold. We
start from the formulation of the randomized control problem
introduced in section \ref{randomizedformulation}. For simplicity of notation,
from now on we drop all superscripts $\hat{}\;$  and start from a setting,
denoted
$(  \Omega,   \calf,  \P,   W,   \mu)$,
where
$(\Omega,   \calf,  \P)$
is a complete probability space,
$  W$
is a standard Wiener process in $\R^{ d}$,
$  \mu=\sum_{n\ge 1}\delta_{(\sigma_n,\eta_n)}$
is a Poisson random measure on  $A$ with intensity $\lambda$, independent of $W$.
We consider  the piecewise constant process $I$ in $A$ associated with $\mu$ defined in
 \eqref{I}, the corresponding trajectory $X$ solution to equation
\eqref{dynXrandom} and the $\P$-complete right-continuous
filtration
$\F^{ W,  \mu}=(\calf^{  W,  \mu}_t)_{t\ge 0}$
generated by  $  W,  \mu$ and defined  by formula
\eqref{filtraznaturaleWmu}.
We recall the estimate
$      \E\,\Big[\sup_{t\in [0,T]}|  X_t|^p\Big]   <\infty
$ for all $p\in [1,\infty)$
(compare
\eqref{EstimateX}).

Our aim is to show that the value of the randomized problem
can be represented in terms of a constrained BSDE, that
we will call \emph{randomized}. From Theorem
\ref{MainThm}
it follows that 
the randomized BSDE 
also represents the value of the original switching problem.

On the space $( \Omega,  \calf, \P)$ equipped with the filtration $\F^{W,\mu}$, let us consider the following constrained BSDE
on the time interval $[0,T]$:
\begin{equation}\label{BSDEconstrained}
\begin{cases}
\vspace{2mm} \dis Y_t \ = \ g( X,  I_T)
+ \int_t^T f_s(X,I_s)\, ds + K_T - K_t - \int_t^TZ_s\,dW_s
- \int_{(t,T]}\!\int_A U_s(a)\,\mu(ds\,da), \\
\dis U_t(a) \ \le \ c_t(X, I_{t-},a).
\end{cases}
\end{equation}

We look for a (minimal) solution to \eqref{BSDEconstrained}  in the sense of the following definition.

\begin{Definition}\label{BSDEdef}
A quadruple $(Y_t,Z_t,U_t(a),K_t)$ $($$t\in [0,T]$, $a\in A$$)$
is called a solution to the BSDE  \eqref{BSDEconstrained} if
\begin{enumerate}
  \item $Y$ $\in$  $\Sc^2(\F^{W,\mu})$, the set of real-valued c\`adl\`ag  $\F^{W,\mu}$-adapted processes  satisfying  $\|Y\|_{\Sc^2}^2$ $:=$ $ \E[\sup_{t\in[0,T]}|Y_t|^2]$ $<$ $\infty$;
  \item $Z$ $\in$ $L_W^2(\F^{W,\mu})$, the set of     $\F^{W,\mu}$-predictable  processes
  with values in $\R^d$
  satisfying $\|Z\|_{L_W^2}^2$ $:=$
      $\E\big[\int_0^T|Z_t|^2dt\big]<\infty$;
  \item $U$ $\in$ $L_{\mu}^2(\F^{W,\mu})$, the set of real-valued $\calp(\F^{W,\mu})\otimes \calb(A)$-measurable processes  satisfying $\|U\|^2_{L_{\mu}^2}$ $:=$
  $\E\big[\int_0^T\int_A|U_t(a)|^2\lambda(da)dt\big]$ $<$ $\infty$;
  \item $K$ $\in$ $\Kc^2(\F^{W,\mu})$, the subset of $\Sc^2(\F^{W,\mu})$ consisting of  $\F^{W,\mu}$-predictable nondecreasing processes with $K_0=0$;
  \item $\P$-a.s. the equality in \eqref{BSDEconstrained}
holds for every $t\in[0,T]$,
  and the constraint $U_t(a)\le c_t(X, I_{t-},a)$ is understood to hold
  $\P(d\omega)\lambda(da)dt$-almost everywhere.
\end{enumerate}
A minimal solution $(Y,Z,U,K)$ is a solution to \eqref{BSDEconstrained} such that for any other solution $(Y',Z'$, $U',K')$, we have
$\P$-a.s., $Y_t\le Y'_t$ for all $t\in [0,T]$.
\end{Definition}

We now state the main result of this section.

\begin{Theorem}
\label{Thm:RandomizedFormula}
There exists a unique minimal solution $(Y,Z,U,K)$ $\in$ $\Sc^2(\F^{W,\mu})\times L_W^2(\F^{W,\mu})\times L_{\mu}^2(\F^{W,\mu})\times\Kc^2(\F^{W,\mu})$ to the randomized BSDE \eqref{BSDEconstrained}. Moreover, we have $Y_0=\sup_{\nu\in\calv} J^\Rc(\nu)$, and, more generally (setting $\eta_0=\xi_0$ for convenience)
\beq
\label{RandomizedFormula}
Y_t &=& \esssup_{\nu\in\calv} \E^\nu\bigg[ \int_t^Tf_s(X,I_s)\,ds +  g(X,I_T)
 -\sum_{n\ge 1}1_{t<\sigma_n<T}
\,c_{ \sigma_n}(  X,  \eta_{n-1}, \eta_n)
 \,\bigg|\, \calf_t^{W,\mu}\bigg].
\enq
\end{Theorem}
\begin{Remark}
{\rm
From Theorems \ref{MainThm} and \ref{Thm:RandomizedFormula} we deduce the BSDE representation for the original optimal switching problem:
\[
Y_0 \ = \ \sup_{\alpha\in\Ac} J(\alpha).
\]
\ep
}
\end{Remark}

\noindent We need the following preliminary result.

\begin{Lemma}
\label{rewritingofcosts}
  For every $\nu\in\calv$ and $t\in [0,T]$, we have $\P$-a.s.
  \beq\label{rewritingofcostsuno}
 \E^\nu\bigg[ \sum_{n\ge 1}1_{t<\sigma_n<T}
\,c_{ \sigma_n}(  X,  \eta_{n-1}, \eta_n)
 \,\bigg|\, \calf_t^{W,\mu}\bigg]&\!=\!&
 \E^\nu\bigg[ \sum_{n\ge 1}1_{t<\sigma_n\le T}
\,c_{ \sigma_n}(  X,  \eta_{n-1}, \eta_n)
 \,\bigg|\, \calf_t^{W,\mu}\bigg]
 \\&\!=\!&\label{rewritingofcostsdue}
 \E^\nu\bigg[
 \int_t^T\int_A c_{ s}(  X,  I_{s-},a)\,\nu_s(a)\;\lambda(da)\,ds
 \,\bigg|\, \calf_t^{W,\mu}
 \bigg].
\enq
In particular, for $t=0$, we have  $J_2^\Rc(\nu)=
 \E^\nu \Big[\int_0^T\int_A c_{ s}(  X,  I_{s-},a)\,\nu_s(a)\;\lambda(da)\,ds\Big]$.
\end{Lemma}

\textbf{Proof.}
The equality in
\eqref{rewritingofcostsuno} is obvious since $\P(\sigma_n=T$ for some $n)=0$ (we note that
in the original switching problem the condition
$\tau_n\neq T$ a.s. had to be imposed: compare
condition (v) in the definition of admissible
switching strategy). Since the $\P^\nu$-compensator of $\mu(ds\,da)$ is $\nu_s(a)\,\lambda(da)\,ds$,
and by {\bf (A1)}-(ii) the random field $c_{ s}(  X,  I_{s-},a)$ is
$\calp(\F^{W,\mu})$ $\otimes$ $ \calb(A)$-measurable and nonnegative,
we obtain
the second equality \eqref{rewritingofcostsdue}:
  \beqs
 \E^\nu\bigg[ \sum_{n\ge 1}1_{t<\sigma_n\le T}
\,c_{ \sigma_n}(  X,  \eta_{n-1}, \eta_n)
 \,\bigg|\, \calf_t^{W,\mu}\bigg]&=&
 \E^\nu\bigg[
 \int_{(t,T]}\int_A c_{ s}(  X,  I_{s-},a) \;\mu(da\,ds)
 \,\bigg|\, \calf_t^{W,\mu}
 \bigg].
 \\&=&
 \E^\nu\bigg[
 \int_t^T\int_A c_{ s}(  X,  I_{s-},a)\,\nu_s(a)\;\lambda(da)\,ds
 \,\bigg|\, \calf_t^{W,\mu}
 \bigg].
\enqs
 \ep

\begin{Remark}
{\rm
It follows from the Lemma that formula  \eqref{RandomizedFormula}
can written
\beqs
Y_t &=& \esssup_{\nu\in\calv} \E^\nu\bigg[ \int_t^Tf_s(X,I_s)\,ds +  g(X,I_T)
 -\int_t^T\int_A c_{ s}(  X,  I_{s-},a)\;\lambda(da)\,ds
 \,\bigg|\, \calf_t^{W,\mu}\bigg].
\enqs
Similar remarks apply to several formulae that follow below.
}
\end{Remark}

\textbf{Proof (of Theorem \ref{Thm:RandomizedFormula})}
Let us introduce for every $n\in\N$ the following penalized BSDE on $[0,T]$:
\begin{equation}\label{BSDEpenalized}
Y_t^n \ = \ g(X,I_T)  + \int_t^T f_s(X,I_s) \,ds + K_T^n - K_t^n
- \int_t^TZ_s^n\,dW_s - \int_t^T\int_A U_s^n(a)\,\mu(ds\,da),
\end{equation}
where
\[
K_t^n \ = \ n \int_0^t \int_A \big(U_s^n(a)-c_s(X,I_{s-},a)\big)^+ \,\lambda(da)\,ds.
\]
    By \eqref{PolGrowth_f_g} and
\eqref{EstimateX} we have
$
\E|g(X,I_T)|^2 < \infty$ and $\E \int_0^T | f_t(X,I_t)|^2\,dt < \infty$, so it follows
 from Lemma 2.4 in \cite{tang_li} that, for every $n\in\N$, there exists a unique solution $(Y^n,Z^n,U^n)$ $\in$ $\Sc^2(\F^{W,\mu})\times L_W^2(\F^{W,\mu})\times L_{\mu}^2(\F^{W,\mu})$ to the above penalized BSDE.

\bigskip

Next we claim that for every $t\in[0,T]$ we have, $\P$-a.s.
\beq
\label{RandomizedFormula_n}
Y_t^n &=& \esssup_{\nu\in\calv_n}
\E^\nu\bigg[ \int_t^T f_s(X,I_s)\,ds + g(X,I_T)-
\sum_{n\ge 1}1_{t<\sigma_n<T}
\,c_{ \sigma_n}(  X,  \eta_{n-1}, \eta_n)
 \,\bigg|\, \calf_t^{W,\mu}\bigg],
\enq
where $\Vc_n=\{\nu\in\Vc\colon\nu\text{ takes values in }(0,n]\}$.
To prove the claim we take any $\nu\in \Vc_n$ and we first notice that
$$
\E^\nu\bigg[ \int_{(t,T]}\int_A U_s^n(a)\,\mu(ds\,da)\,\bigg|\, \calf_t^{W,\mu}\bigg]=
\E^\nu\bigg[ \int_t^T\int_A U_s^n(a)\,\nu_s(a)\,\lambda(da)\,ds\,\bigg|\, \calf_t^{W,\mu}\bigg],
$$
because the $\P^\nu$-compensator of $\mu(ds\,da)$ is $\nu_s(a)\lambda(da)\,ds$.
Next we note that the process  $\int_0^\cdot Z^n_s\,dW_s $ is a $\P^\nu$-local
martingale, since $W$ is a Wiener process under $\P^\nu$;
recalling that $d\P^\nu=\kappa^\nu_Td\P$ and using the estimates
\eqref{stimasukappa} it is easy to prove that it is in fact a $\P^\nu$-martingale,
so that in particular
$$
\E^\nu\bigg[ \int_t^TZ^n_s\,dW_s \,\bigg|\, \calf_t^{W,\mu}\bigg]=0.
$$
So taking expectation   $\E^\nu$ in \eqref{BSDEpenalized},
adding and subtracting both sides of equality \eqref{rewritingofcostsdue}
and rearranging terms we obtain
\beq
\label{relazionefondamentale}
Y_t^n &=&   \E^\nu\bigg[ \int_t^Tf_s(X,I_s)\,ds +  g(X,I_T)
 -\sum_{n\ge 1}1_{t<\sigma_n<T}
\,c_{ \sigma_n}(  X,  \eta_{n-1}, \eta_n)
 \,\bigg|\, \calf_t^{W,\mu}\bigg]
 \\\nonumber
 &+&
 \E^\nu\bigg[
     \int_t^T \int_A \Big\{n\big(U_s^n(a)-c_s(X,I_{s-},a)\big)^+ -
     \big(U_s^n(a)-c_s(X,I_{s-},a)\big)\nu_s(a)
     \Big\}\,\lambda(da)\,ds
 \,\bigg|\, \calf_t^{W,\mu}\bigg].
\enq
This is sometimes called the fundamental relation for the penalized control
problem corresponding to admissible controls $\Vc_n$.
The term in curly brackets
$\Big\{\ldots
     \Big\}$ is nonnegative, since $\nu_s(a)$ takes values in $(0,n]$ and
     we have the numerical inequality $nx^+\ge x\nu$ for every $x\in\R$
     and $\nu\in (0,n]$.  It follows that
\beq
\label{relazfondineq}
Y_t^n &\ge &
\E^\nu\bigg[ \int_t^T f_s(X,I_s)\,ds + g(X,I_T)-
\sum_{n\ge 1}1_{t<\sigma_n<T}
\,c_{ \sigma_n}(  X,  \eta_{n-1}, \eta_n)
 \,\bigg|\, \calf_t^{W,\mu}\bigg], \quad \nu\in\Vc_n.
\enq
Now we show that the term in curly brackets can be made as small as we wish for
an appropriate choice of $\nu\in\Vc_n$.      We note that,
given $0<\eps <n$ and $x\in\R$ and choosing
$$
 \bar\nu=n\,1_{\{x\geq0\}} + \eps\,1_{\{-1<x<0\}}
 - (\eps/x) \,1_{\{x\leq -1\}}$$
 we have $\bar\nu\in [\eps,n]$ and $nx^+-x\bar\nu\le \eps$.
     So it follows that     setting
 \beqs\nu_s^{\eps,n}(a)&=&n\,1_{\{U_s^n(a) -c_s(X,I_{s-},a)\geq0\}}
 + \eps\,1_{\{-1<U_s^n(a) -c_s(X,I_{s-},a)<0\}}
 \\&&
  - \eps\,(U_s^n(a) -c_s(X,I_{s-},a))^{-1}\,1_{\{U_s^n(a) -c_s(X,I_{s-},a)\leq -1\}},
 \enqs
 we have $\nu^{\eps,n}\in \Vc_n$ and
    \[
\Big\{n\big(U_s^n(a) -c_s(X,I_{s-},a)\big)^+
- \nu_r^{\eps,n}(a)( U_s^n(a)-c_s(X,I_{s-},a)) \Big\}\leq  \eps
\]
($\nu_s^{\eps,n}(a)$ is an approximation
of $n\,1_{\{U_s^n(a) -c_s(X,I_{s-},a)\geq0\}}$ which is not in $ \Vc_n$
since it can take the value zero). From
\eqref{relazionefondamentale} it follows that
\beqs
Y_t^n &\le &
\E^{\nu^{\eps,n}}\bigg[ \int_t^T f_s(X,I_s)\,ds + g(X,I_T)-
\sum_{n\ge 1}1_{t<\sigma_n<T}
\,c_{ \sigma_n}(  X,  \eta_{n-1}, \eta_n)
 \,\bigg|\, \calf_t^{W,\mu}\bigg]+\eps\,(T-t)\,\lambda(A),
\enqs
which, together with \eqref{relazfondineq},  proves
the claim
   \eqref{RandomizedFormula_n}.

\bigskip

Recalling that $d\P^\nu=\kappa^\nu_Td\P$, using the estimates
\eqref{stimasukappa}, \eqref{stimasureward} and recalling
 \eqref{PolGrowth_f_g} and \eqref{EstimateX}, we deduce that
\begin{equation}
\label{YnUpperBound}
\sup_n Y_t^n \ < \ \infty, \qquad \text{for all }0\leq t\leq T.
\end{equation}

Let us define
  $\check{g}$, $ \check{Y}$   and $\check{U}$ by the equalities
 \[ \left\{ \begin{array}{l}
 \dis  \check{Y}_t = Y_t - \int_{0}^t\int_A  c_s(X,I_{s-},a)\, \mu(ds,da), 
 \\
  \dis \check{g} = g(X,I_T) - \int_{0}^T\int_A c_s(X,I_{s-},a)\, \mu(ds,da) \\
 \dis \check{U}_t(a)  = U_t(a) - c_t(X,I_{t-},a),
\end{array} \right.
 \]
so that equation \eqref{BSDEconstrained} can be written as follows:  
\begin{equation}\label{BSDEconstrainedbis}
\begin{cases}
\vspace{2mm} \dis \check{Y}_t  = \check{ g}+ \int_t^T \bigg( f_s(X,I_s) - \int_A \check{U}_s(a)\,\lambda(da)\bigg) ds + K_T - K_t 
- \int_t^T Z_s\,dW_s - \int_t^T \int_A \check{U}_s(a)\,\tilde\mu(ds\,da), \\
\dis \check{U}_t(a) \ \le \ 0,
\end{cases}
\end{equation}
where $\tilde \mu(dt\,da) = \mu(dt\,da) - \lambda(da)\,dt$ denotes the compensated
Poisson measure.
Let us check that 
$(Y,Z,U,K)$ belongs to the space $\Sc^2(\F^{W,\mu})\times L_W^2(\F^{W,\mu})\times L_{\mu}^2(\F^{W,\mu})\times\Kc^2(\F^{W,\mu})$
if and only if $(\check{Y},Z,\check{U},K)$ does. In fact, noting that
the process 
$ c_t(X,I_{t-},a)$  is $\calp(\F^{W,\mu})\otimes \calb(A)$-measurable
and non-negative, it is enough to verify that
$$
\E\left[\left| \int_{0}^T\int_A c_s(X,I_{s-},a)\, \mu(ds,da)
\right|^2\right]<\infty, $$
and
$$\E\left[ \int_{0}^T\int_A c_s(X,I_{s-},a)^2\, \lambda(da)\,ds\right]<\infty.
$$
These inequalities follow from the growth assumption \eqref{PolGrowth_f_g}
and the estimate \eqref{EstimateX}, taking into account that the
random measure $\mu((0,T]\times A)$ has Poisson law with
parameter $\lambda(A)T$. It also follows that $\check{g}$ also belongs to $L^2$ and it is
$\F^{W,\mu}$-measurable.
We conclude that $(Y,Z,U,K)$ is the minimal solution to \eqref{BSDEconstrained}
if and only if $(\check{Y},Z,\check{U},K)$
is the minimal solution to \eqref{BSDEconstrainedbis}. 

Next we also note that equation \eqref{BSDEconstrainedbis}
is a particular case of a backward stochastic differential equation studied in a 
general non-Markovian framework in \cite{KP12}. In particular, 
existence and uniqueness of the minimal solution to equation \eqref{BSDEconstrainedbis}
(or, equivalently, to equation \eqref{BSDEconstrained}) follow from Theorem 2.1 in \cite{KP12}.
Indeed, Assumption {\bf (H0)} in \cite{KP12} is clearly satisfied. Concerning Assumption {\bf (H1)}, this is only used in Lemma 2.2 of \cite{KP12}
to prove that the sequence $(Y^n)_n$ satisfies \eqref{YnUpperBound},
a property that in our setting has been proved by different arguments.
Finally, from Theorem 2.1 in \cite{KP12} we also have that $Y_t^n(\omega)$ converges increasingly to $Y_t(\omega)$ as $n\rightarrow\infty$,
$\P(d\omega)$-a.s.
Since $\Vc=\cup_n\Vc_n$, letting $n\rightarrow\infty$ in \eqref{RandomizedFormula_n} we obtain \eqref{RandomizedFormula}.
\qed

\begin{Remark}
{\rm
It follows from the previous proof that the process
$\nu^{\eps,n}$ constructed above  satisfies
$J^\Rc(\nu^{\eps,n})+\delta >Y_0=\sup_{\nu\in\calv} J^\Rc(\nu)$
for arbitrary $\delta>0$, provided $\epsilon$ is sufficiently
small and $n$ sufficiently large: in other words, it is
$\delta$-optimal for the randomized problem. One
could then repeat the arguments of section
\ref{firstineq} and construct a corresponding control (call it  
$ \alpha^{\delta}$)
which is $\delta$-optimal for the non-randomized 
switching
problem. 
However, this is not enough to conclude that 
$ \alpha^{\delta}$ is $\delta$-optimal for the
original switching problem: indeed, as specified in 
Proposition
\ref{P:Ineq_I}, it only belongs to the class 
$ \hat\Ac^{\infty}$ introduced
in \eqref{primalvalueextinfty} and so it is not adapted
to the Brownian filtration of the original problem.
}
\end{Remark}

\vspace{3mm}

Formula \eqref{RandomizedFormula} shows that the process $Y$
constructed in Theorem \ref{Thm:RandomizedFormula}
can be seen as the value
of an optimization problem. Our final result shows that it satisfies
a version of dynamic programming principle (DPP in short) in the randomized context.
We omit the proof which is very similar to Lemma 4.8 in \cite{FP15}
or Theorem 5.3 of \cite{BCFP16bAAP},
after obvious changes of notation.

\begin{Theorem}\label{T:DPPRandomizedProblem}
For all $0\leq t\leq T$, we have
\beq\label{DPPRandomizedProblem}
Y_t &=& \esssup_{\nu\in\calv} \esssup_{\tau\in\Tc_t} \E^\nu\bigg[
\int_t^\tau f_r(X,I_r)\,dr
-
\sum_{n\ge 1}1_{t<\sigma_n<\tau }
\,c_{ \sigma_n}(  X,  \eta_{n-1}, \eta_n)
+ Y_\tau \,\bigg|\, \calf_t^{W,\mu}\bigg] \notag \\
&=& \esssup_{\nu\in\calv} \essinf_{\tau\in\Tc_t} \E^\nu\bigg[
\int_t^\tau f_r(X,I_r)\,dr
-
\sum_{n\ge 1}1_{t<\sigma_n< \tau}
\,c_{ \sigma_n}(  X,  \eta_{n-1}, \eta_n)
+ Y_\tau \,\bigg|\, \calf_t^{W,\mu}\bigg],
\enq
where $\mathcal T_t$ denotes the class of $[t,T]$-valued $\F^{W,\mu}$-stopping times.
\end{Theorem}

\medskip
\noindent {\bf Remark 5.4.}
Besides its intrisinc interest, this result may also be used
in the Markovian context (i.e., when the coefficients
are not path-dependent) in order to give an alternative
proof that the value function (of the primal optimal switching problem) is a viscosity solution to
the corresponding Hamilton-Jacobi-Bellman equation (1.3). In our present framework this may lead to extending
known results on (1.3) to the case
of infinitely many modes.\\ 
\noindent However, such a randomized DPP is not \textit{a priori} equivalent to the standard one: more precisely, the stopping times involved in this randomized DPP have to be 
$\mathbb{F}^{W, \mu}$-adapted. Thus, similarly as in Remark 5.3, they are not necessarily adapted to the original brownian filtration. Additionnally, it does not seem clear that all ingredients already used in the proof of Theorem 5.1 in [11] can be reproduced in the present setting. In the present case, one has to deal with the sum of penalty costs which plays a prominent role in particular in the constraint condition. Indeed, the jump constraint (appearing in the so-called \textit{randomized} BSDE) is naturally related with the obstacle in the (infinite) PDE system and one can guess that such an issue has to be dealt with. For this reason we left it to future study.

\paragraph*{Acknowledgments}
$\\$
The two authors thanks both Le Mans Universit\'e (France) and Universit\`a degli Studi Di Milano (Italy) for their kind hospitality during their respective visits. We also adress special thanks to Prof. Sa\"id Hamad\`ene (Le Mans Universit\'e, Le Mans) and Prof. Idris Kharroubi (Sorbonne Universit\'e, Paris) for fruitful discussions on related topics and the two anonymous referees for their comments which helped us to improve this paper.


\end{document}